\documentclass[11pt, leqno,amscd, psamsfonts]{amsart}
\usepackage{amsthm, amsmath}
\usepackage{amssymb} 
\usepackage{amscd}
\def\arrowdown#1#2{\Big\downarrow \rlap{$\vcenter{\hbox{$\scriptstyle#2$}}$}
{\hbox to -10pt{\hss{$\vcenter{\hbox{$\scriptstyle#1$}}$}}}}
\def\arrowup#1#2{\Big\uparrow \rlap{$\vcenter{\hbox{$\scriptstyle#2$}}$}
{\hbox to -10pt{\hss{$\vcenter{\hbox{$\scriptstyle#1$}}$}}}}
\DeclareFontEncoding{OT2}{}{} % to enable usage of cyrillic fonts
  \newcommand{\textcyr}[1]{%
    {\fontencoding{OT2}\fontfamily{wncyr}\fontseries{m}\fontshape{n}%
     \selectfont #1}}
\newcommand{\Sha}{{\mbox{\textcyr{Sh}}}}

\newtheorem{thm}{Theorem}[section]
\newtheorem{prop}[thm]{Proposition}
\newtheorem{lemma}[thm]{Lemma}

\newtheorem{Definition}[thm]{Definition}
\newtheorem{Remarknumb}[thm]{Remark}

\newtheorem{Example}[thm]{Example}

\newtheorem{conjecture}[thm]{Conjecture}
\newtheorem{locconj}[thm]{Localization Conjecture}
\newtheorem{cor}[thm]{Corollary}

\newcounter{ex}[section]

\numberwithin{equation}{section}
\numberwithin{thm}{section}

\newcommand{\cal}{\mathcal}

\newcommand{\E}{{\mathcal E}}
\newcommand{\A}{{\mathcal A}}

\newcommand{\Cl}{{\rm Cl}}

\newcommand{\C}{{\bf C}}
\newcommand{\R}{{\rm R}}

\newcommand{\Q}{{\bf Q}}

\newcommand{\Hom}{{\rm Hom}}
\newcommand{\X}{{\mathcal X}}

\newcommand{\Gg}{{\mathcal G}}
\newcommand{\Hh}{{\mathcal H}}
\newcommand{\Gm}{{{\bf G}_m}}
\newcommand{\Z}{{\bf Z}}
\newcommand{\F}{{\mathcal F}}

\newcommand{\ti}{\tilde}
\newcommand{\Spec}{{\rm Spec }\, }
 \renewcommand{\O}{{\mathcal O}}

\newcommand{\M}{{\mathcal M}}
\renewcommand{\L}{{\mathcal L}}

\newcommand{\Gr}{{\rm G}}

\newcommand{\Kr}{{\rm K}}

\newcommand{\Tr}{{\rm Tr}}

 \newcommand{\fp}{{{\bf F}_p}}

\newcommand{\Pic}{{\rm Pic}}
\newcommand{\PP}{{\mathcal P}}

\newcommand{\acfp}{\overline{{\bf F}}_p}

\renewcommand{\mod}{\,\hbox{\rm mod}\,}

\def\thfill{\null\nobreak\hfill}

\def\endproof{\thfill\vbox{\hrule
  \hbox{\vrule\hbox to 5pt{\vbox to 5pt{\vfil}\hfil}\vrule}\hrule}}

\renewcommand{\P}{{\cal P}}
\newcommand{\Hr}{{\rm H}}

\begin{document}
\title[Cubic structures]{Cubic structures, equivariant Euler characteristics and 
lattices of modular forms}

\author[T. Chinburg]{T. Chinburg*}\thanks{*Supported by NSF Grant \# DMS00-70433 and \# DMS05-00106}
\address{Ted Chinburg, Dept. of Math\\Univ. of Penn.\\Phila. PA. 19104, U.S.A.}
\email{ted@math.upenn.edu}

\author[G. Pappas]{G. Pappas\dag}
\thanks{\dag Supported by NSF Grants \# DMS05-01049 and \# DMS01-11298 (via the Institute for Advanced Study)}
\address{Georgios Pappas, Dept. of Math\\Michigan State Univ.\\East Lansing, MI 48824, U.S.A.}
\email{pappas@math.msu.edu}

\author[M. Taylor]{M. J. Taylor\ddag}
\thanks{\ddag EPSRC Senior Research Fellow and Royal Society Wolfson Merit award}
\address{Martin J. Taylor, Dept. of Math\\University of Manchester, M60 1QD, U.K.}
\email{Martin.Taylor@manchester.ac.uk}

\date{\today}

\begin{abstract}
We use the theory of cubic structures to give a fixed point
Riemann-Roch formula for the equivariant Euler characteristics of
coherent sheaves on projective flat schemes over $\bf Z$ with a 
tame action of a finite abelian group. This formula supports a 
conjecture concerning the extent to
which such equivariant Euler characteristics  may be determined from the
restriction of the sheaf to an infinitesimal neighborhood of the fixed
point locus.  Our results are applied to study the module structure of
modular forms having Fourier coefficients in a ring
of algebraic integers, as well as the action of diamond Hecke
operators on the Mordell-Weil groups and Tate-Shafarevich groups of
Jacobians of modular curves.
\end{abstract}

\maketitle

\centerline{\sc Contents}
\medskip

\noindent \S 1. Introduction\\
\S 2. Cubic structures and categories\\
\S 3. The main calculation and result\\
\S 4. Galois structure of modular forms\\
\S 5. An equivariant Birch and Swinnerton-Dyer  relation
\medskip

%\tableofcontents
\section{Introduction}
\label{s:intro}

Let $X$ be a projective flat scheme over 
$\Spec(\Z)$ having an action of a finite group $G$.
Let $\F$ be a $G$-equivariant coherent locally free sheaf on $X$.
One then has an Euler characteristic
$$
\chi(X,\F)=\sum_i(-1)^i[{\rm H}^i(X, \F)]
$$
in the Grothendieck group ${\rm G}_0({\bf Z}[G])$ of all finitely 
generated $G$-modules.  If the action of $G$ on $X$ is tame,
as we shall assume for most of the article,
there is a refinement $\chi^P(X,\F)$ of $\chi(X,\F)$ in the Grothendieck
group  ${\rm K}_0({\bf Z}[G])$ of all finitely generated
projective ${\bf Z}[G]$-modules.  Let $X^g$ be the subscheme of $X$ fixed by the action of $g\in G$,
 and let $X' = \cup_{  e\ne g \in G} X^g$.  The goal of this paper is
to compute $\chi^P(X,\F)$ and $\chi(X,\F)$ when $G$ is abelian via the restriction 
of $\F$ to an infinitesimal neighborhood of $X'$.  In view of the fact that in general,
${\rm G}_0({\bf Z}[G])$ and ${\rm K}_0({\bf Z}[G])$ have non-trivial
torsion subgroups,
this can be viewed
as a variant of the problem of finding a Lefschetz-Riemann-Roch
formula for  $\chi^P(X,\F)$ and $\chi(X,\F)$ with an explicit, bounded denominator.  
Some applications we will discuss
are to the Galois module structure of lattices of modular forms, and of the Tate-Shafarevich
groups and Mordell-Weil groups of Jacobians of modular curves, with respect to the
action of diamond Hecke operators.

Our main tool is the theory of cubic structures. These
were first studied in detail by 
Breen in [Br], building on work of  Mumford and Grothendieck concerning biextensions
and the Theorem of the Cube for abelian varieties. The use of cubic structures in 
our context was first introduced by
one of us in [P1], [P2] and [P3] to prove refined Galois module structure results
when the action of the group $G$ on $X$ is free, that is when $X'= \emptyset$. 
Here we present a reformulation 
of cubic structures which is tailored to proving Lefschetz-Riemann-Roch results when $G$ is abelian
and the action of $G$ on $X$ is tame.  Our results improve on what can be shown
using the localization techniques of Thomason [T1] and [T2], which generally
determine Euler characteristics in   ${\rm G}_0({\bf Z}[G])$  or ${\rm K}_0({\bf Z}[G])$ 
only up to classes of finite order dividing an unspecified power
of the order of $G$. (Thomason actually deals with the case in which 
the action is by a multiplicative group scheme, but some of his techniques
extend to the case of constant groups $G$; see [CEPT3].) 

We will study the image $\overline{\chi}^P(X,\F)$ of $\chi^P(X,\F)$
in the classgroup $\Cl(\Z[G])$ of $\Z[G]$, which is the quotient ${\rm K}_0(\Z[G])$ by the subgroup generated by the class of $\Z[G]$.  The projection ${\rm K}_0(\Z[G]) \to \Cl(\Z[G])$
identifies the torsion subgroup of ${\rm K}_0(\Z[G])$ with the finite abelian group $\Cl(\Z[G])$.  With this identification, $\chi^P(X,\F) = \overline{\chi}^P(X,\F) + r(\Z[G])$ where $$  r =  \frac{1}{\#G} \ {\rm rank}_\Z\, ( \chi^P(X,\F)) =   \frac{1}{\#G}  \ \sum_i(-1)^i\ {\rm dim}_{\Q}({\rm H}^i(X_\Q,\F_\Q))$$ is an integer which can be determined using a non-equivariant Riemann Roch formula
 for the Euler characteristic of the restriction $\F_\Q$ of $\F$ to the general fiber $X_\Q$ of $X$.

The results we obtain, along with those of [P3],
support the following conjecture about $\overline{\chi}^P(X,\F)$.
Define $\hat X'$ to be the formal completion of $X'$ in $X$, and
let $\F{|\hat X'}$ be the $G$-sheaf on $\hat X'$ which is the restriction of $\F$. 
Since $X$ is projective and the action of $G$ on $X$ is tame, the image
of $X'$ in $\Spec(\Z)$ is  a proper closed set.
\smallskip

\begin{locconj}
\label{con:localconj}
There is an integer $N \ge 1$ which depends only on ${\rm dim}(X)$, and an integer 
$\delta \ge 1$ which depends only on $\# G$, for which the
following is true.    Let
$m = {\rm g.c.d.}(N,\#G)^{\delta}$.  
\smallbreak
\begin{enumerate}
\item[a.] {\bf (Input localization)} The class $m \cdot \overline{\chi}^P(X,\F)$ depends only
on the pair $(\hat X', \F{|\hat X'})$.
\medbreak
\item[b.]{\bf (Output localization)}  The class $m\cdot  \overline{\chi}^P(X,\F)$ is the class of a projective ideal $I \subset \Z[G]$ such that the order of $\Z[G]/I$ is supported on the image of $X'$ in $\Spec(\Z)$. 
\end{enumerate}
\end{locconj}
\smallskip

\noindent The integer $m$ can be thought of as the denominator in a Riemann-Roch
formula for $\overline{\chi}^P(X,\F)$. 
We have referred to the two parts of this conjecture in this way because 
 part (a) restricts the amount of input information needed to
determine $m \cdot \overline{\chi}^P(X,\F)$, while part (b) restricts the 
complexity of the output information
needed to express  $m\cdot  \overline{\chi}^P(X,\F)$.  

We would like to thank L. Illusie for
pointing out a parallel between Conjecture \ref{con:localconj} and results of Deligne, Abb\`es and Vidal
concerning the equality of $l$-adic \'etale Euler characteristics of constructible sheaves
having the same rank and wild ramification at infinity.  For the precise statement of these
results, see  [I], [A] and [V].

The  methods of [T1], [T2] and [CEPT3] lead to a proof of Conjecture \ref{con:localconj}(a) modulo
classes annihilated by a power of $\# G$.  Our main results imply the following theorem.
\begin{thm}
\label{thm:inputoutput}
Suppose $G$ is abelian and that the branch locus of $X \to X/G = Y$
is supported off the prime divisors of $\# G$. Suppose also that $Y$ is regular and $X$ is normal.
\begin{enumerate}
\item[a.]  The input localization conjecture is true for $X$.  If $dim(X) \le 4$,
one can let  $m = 1$ if $\# G$ is odd, and $m = 2$
if $\# G$ is even.
\item[b.] If $m$ is the integer used in part (a) above, then the class  $2|G|^{dim(X)}m\cdot  \overline{\chi}^P(X,\F)$ 
is the class of a projective ideal $I \subset \Z[G]$ such that the order of $\Z[G]/I$ 
is supported on the image of $X'$ in $\Spec(\Z)$. 
\end{enumerate}
\end{thm}
We will in fact show that for abelian $G$, the output localization conjecture holds for the image
of $\chi^P(X,\F)$ under a map $\Theta_n$ described below. 

We now
explain the appearance of the integers $N$ and $\delta$ in Conjecture \ref{con:localconj}.
Suppose the action of $G$ on $X$ is free, so that  $X'=\emptyset$.
In this case, each part of the conjecture is separately equivalent to $m \cdot \overline{\chi}^P(X,\F) = 0$.  
This statement for an explicit $N$ and $\delta$ depending only on $\dim(X)$ and $\#G$, 
respectively, was proved in [P3] when all the Sylow subgroups of $G$ are abelian.
However, in [P3] it was suggested that
Vandiver's conjecture (which claims that the class number of $ \Z[\zeta_p+\zeta_p^{-1}]$
is prime to $p$) fails for exactly those primes $p$
for which there is an $X$ with free $\Z/p\Z$-action and $p>{\rm dim}(X)$ 
with $\overline{\chi}^P(X,\F)\ne 0$.  We have introduced the integer $N$ in 
the Conjecture because if one thought the conjecture were true
with $N = 1$, then [P3] suggests that this would imply the truth
of Vandiver's conjecture, and this seems too strong.   

Theorem \ref{thm:inputoutput} is a consequence of the explicit formula 
in Theorem \ref{mainthm}, which is the main result of this paper.  We now
indicate some consequences of this formula concerning modular forms,
Mordell-Weil and Tate-Shafarevich groups.  

Let $p\equiv 1 $ mod $24$ be a prime and  $\Gamma = (\Z/p\Z)^*/\{\pm 1\}$.  Suppose $\chi:\Gamma \to \mu_r\subset \Z[\zeta_r]^*$
is a character of prime order $r|(p-1)$ with $r > 3$. 
Let $S_2(\Gamma_1(p),\Z[\zeta_r])_\chi$ be the $\Z[\zeta_r]$-module of cusp forms  
$F(z)=\sum_{n\geq 1}a_ne^{2\pi i n z}$
of weight $2$, level $p$ and Nebentypus character $\chi$ 
whose Fourier coefficients $a_n$  belong to $\Z[\zeta_r]$.  Then 
$S_2(\Gamma_1(p),\Z[\zeta_r])_\chi$ is a locally free $\Z[\zeta_r]$-module,
of rank $n(\chi) = (p - 25)/12$ by the classical Chevalley-Weil
theorem.  (See [CW];  this was one of the first  coherent
Lefschetz-Riemann-Roch theorems.)   For $a \in (\Z/r\Z)^*$ let $\{a\}$ be the unique integer in the range $0 < \{a\} < r$ having residue class $a$, and let $\sigma_a \in {\rm Gal}(\Q(\zeta_r)/\Q)$ be the automorphism for which
$\sigma_a(\zeta_r) = \zeta_r^a$. Define $\omega_r:(\Z/r\Z)^* \to \Z_r^*$ to
be the Teichm\"uller character.  We embed $\Z$ (resp. $\Z_r$) into the
profinite completion $\hat{\Z} = \prod_{l\ {\rm prime}}\Z_l$ of $\Z$ diagonally (resp. 
via the factor $l = r$).  We then have a modified quadratic Stickelberger
element of $\hat{\Z}[{\rm Gal}(\Q(\zeta_r)/\Q)]$ defined by
$$
\theta_2 =  \sum_{a  \in (\Z/r)^*}  \frac{ (p-1)}{24r^2} \left (\{a\}^2 - \omega_r( a)^2\right )\sigma^{-1}_a.
$$
(For the standard definition of the quadratic Stickelberger element see for example [KL, p.115],
also [CNT, p.308].)  Since the ideal class group $\Cl(\Z[\zeta_r])$ 
is finite, $\theta_2$ acts on this group.  Let $\PP_\chi$ be the prime ideal of $\Z[\zeta_r]$ over $(p)$ 
with the property that the reduction of $\chi$ modulo $\P_\chi$ is the
$\frac{p-1}{r}$ power of the identity character $(\Z/p\Z)^*\to {\bf F}_p^*$.  

Let $X_1(p)_\Q$ be the canonical model over $\Q$ of the modular curve, and identify
$\Gamma$ with the group of diamond Hecke operators acting on $X_1(p)_\Q$.  For each subgroup $H \subset \Gamma$ let $X_{H,\Q} = X_1(p)_\Q/H$.  By applying Theorem \ref{mainthm} to a suitable integral model $X_H \to X_0$ of the $G = \Gamma/{\rm ker}(\chi)$ cover $ X_{{\rm ker}(\chi),\Q} \to X_{\Gamma,\Q} = X_0(p)_\Q$, we will prove the following result: 
\begin{thm}
\label{qua1}
Suppose ${\mathfrak A}\subset \Z[\zeta_r]$ is an ideal with ideal class $\theta_2\cdot [\PP_\chi]$. 
Then we have $$S_2(\Gamma_1(p),\Z[\zeta_r])_\chi\simeq \Z[\zeta_r]^{n(\chi) - 1}\oplus {\mathfrak A}$$
as $\Z[\zeta_r]$-modules.
\end{thm}

\noindent This theorem verifies
Conjecture \ref{con:localconj} for $X_H \to X_0$ with $N = 1 = \delta$.

Let $J_H$ be the Jacobian of  $X_{H,\Q}$, and define $J_H(\Q)$ and $\Sha (J_H)$ to be the Mordell-Weil and
Tate-Shafarevich groups of $J_H$ over $\Q$.   
We will assume that  $\Sha (J_H)$  is finite. 
 For $G = \Gamma/{\rm ker}(\chi)$ as above, tensoring $\Z[G]$-modules with the 
ring homomorphism $\Z[G]\to \Z[\zeta_r,\frac{1}{r}]$  induced by $\chi$
gives a Steinitz class homomorphism 
$$
s_\chi:\Gr_0(\Z[G]) \to \Gr_0(\Z[\zeta_r,\frac{1}{2r}])/\{{\rm free \ modules}\} = \Cl(\Z[\zeta_r,\frac{1}{2r}])= \Cl(\Z[\zeta_r,\frac{1}{2}]).
$$

In \S \ref{BSD} we describe an equivariant version of the Birch Swinnerton-Dyer conjecture.
This should follow from the equivariant Tamagawa number
conjecture resulting from the work of Bloch and Kato, Fontaine and Perrin-Riou
and  of Burns and Flach (see [F] and Remark \ref{bigrem}(a)).

\begin{thm}
\label{qua2}
If the Birch and Swinnerton-Dyer conjecture of \S \ref{BSD} holds
then
\begin{equation}\label{bsd3}
 \overline {\theta_2[\PP_\chi]} = s_\chi( [\Sha(J_H)])-\overline{s_\chi( [J_H(\Q)])}-
s_\chi( [J_H(\Q)]) 
\end{equation}
in $\Cl(\Z[\zeta_r, \frac{1}{2}])$, where $\overline{\mathcal{D}}$ is the complex
conjugate of an ideal class $\mathcal{D}$.  
\end{thm}

Now let $C(p)$ be the subgroup of $J_1(p) = J_{\{e\}}$ generated by differences
of $\Q$-rational cusps of $X_1(p)_\Q$;  these are the cusps over the cusp
$\infty$ of $X_0(p)_\Q$. (Conrad, Edixhoven and Stein 
conjecture in [CES] that $C(p) = J_1(p)(\Q)_{\rm tor}$.)
Let $J'_H(\Q)$ be the image of $J_H(\Q)$ in $J_1(p)(\Q)/C(p)$ under the pullback 
homomorphism associated to the quotient morphism $X_1(p)_\Q \to X_{H,\Q}$.
Using a result of Kubert-Lang, we show that $s_\chi([J_H(\Q)])=s_\chi( [J'_H(\Q)])$.
This implies the following.

\begin{cor}
\label{corBSD}
Suppose $ \theta_2 [\PP_\chi] \ne 0 $ in $\Cl(\Z[\zeta_r,\frac{1}{2}])$.  Then at least one of $s_\chi([J'_H(\Q)])$
or $s_{\chi}([\Sha (J_H)])$ is non-trivial if the Birch Swinnerton-Dyer
conjecture of \S \ref{BSD} holds.  Suppose in addition that $C(p) = J_1(p)(\Q)_{\rm tor}$,
as conjectured in [CES].  Then either the $\chi$-eigenspace of
$J_1(p)(\Q)\otimes_\Z{\bf C}$ is non-trivial, or
 $s_\chi([\Sha(J_H)])$ is non-trivial.  
\end{cor}

Whether or
not $ \theta_2 [\PP_\chi]$ is trivial in $\Cl(\Z[\zeta_r,\frac{1}{2}])$ depends only on $p$ and not on the choice of the
character $\chi$ of order $r$.   
If $r$ is fixed and $\theta_2(\Cl(\Z[\zeta_r,\frac{1}{2}])) \ne 0$, e.g. if $r = 191$, the 
set of $p$ for which $ \theta_2 [\PP_\chi] \ne 0$ has a positive Dirichlet density.  
It is not clear, though, whether the alternatives in  Corollary \ref{corBSD} for
such $p$ should occur with a well-defined density, or whether one should occur
asymptotically more often than the other. Nevertheless, one should mention the 
following:
The L-functions of cusp eigenforms of weight $2$, level $p$ and
Nebentypus character $\chi$ are not self-dual. According to the expectations of Katz-Sarnak ([KS]),
the percentage of such forms of level $p\leq x$ whose L-functions
do not vanish at $s=1$ should approach $100\%$ as $x\to\infty$. Hence, it is plausible
that the $\chi$-eigenspace of $J_1(p)(\Q)\otimes{\bf C}$ is trivial for many primes $p$.   
However, since there could be many Galois orbits of such
eigenforms at each level $p$ 
the heuristics of [KS] do not actually imply that this is the case  
for an infinite number of primes.

The principle behind Theorem \ref{qua2} is that Birch Swinnerton-Dyer conjectures
predict identities in $\Gr_0(\Z[\Gamma])$ involving classes constructed
from L-series derivatives, height pairings, period maps, de Rham and Betti cohomology groups,
Mordell-Weil groups and Tate-Shafarevich groups.  Combinations of
these classes which can be defined by Galois-equivariant functions from the
characters of $\Gamma$ to the complex numbers should be  trivial or nearly so in $\Gr_0(\Z [\Gamma])$. 
(This accounts for the fact that values of L-series derivatives, regulators and period maps do not appear in
Theorem \ref{qua2}.) Determining some of the other classes, e.g. the ones associated to de Rham and Betti cohomology over $\Z$,
then leads to predictions about the others, as in Theorem \ref{qua2}.  It would
be of great interest to find a motivic cohomology formalism which would make such predictions
directly, without the intervention of L-value conjectures;
such conjectures play no role in Theorems \ref{mainthm} and \ref{qua1}.

Let us remark here that the  methods used to show Theorem \ref{qua1} apply in many other situations.
We plan to consider generalizations  to other 
spaces of weight two modular forms and covers of modular curves in future work.
An extension to higher weight modular forms is considered by E. Gurel in his  
thesis. Shimura varieties more general than modular curves also provide with many 
examples of Galois covers to which our techniques can be applied.
This way one can hope to obtain information on the Galois module structure 
of lattices of automorphic forms for other groups besides ${\rm GL}(2)/\Q$.
 
We now explain our approach to computing  $\overline{\chi}^P(X,\F)$ for general $X$ assuming that the quotient $Y=X/G$ is a regular scheme which is projective and flat over $\Spec(\Z)$ and
of relative dimension $d$.  Our main tools are the existence of cubic structures on the determinant of cohomology
and a localized Riemann-Roch theorem. For simplicity, we will restrict our discussion 
to the case $\F=\O_X$. The image of $\overline{\chi}^P(X,\O_X)\in \Cl(\Z[G])$ maps, under the natural isomorphism $\Cl(\Z[G]) \xrightarrow{\sim} \Pic(\Z[G])$ given by the determinant of $\Z[G]$-modules,
to the class of the  locally free rank one $\Z[G]$-module given by the determinant of cohomology
$\det \R\Gamma(X,\O_X)$ in $\Pic(\Z[G])$.  

Suppose that $M$ is a locally free rank one 
left $\Z[G]$-module. To each subset $I \subset \{1,\ldots, n\}$
we associate a homomorphism $I:G \to G^{n}$ from $G$ to the product $G^{n}$ 
of $n$ copies of $G$ which sends $g \in G$ to the element of $G^{n}$ having
$g$ in the $i^{th}$ component for $i \in I$ and the identity element of $G$
in all other components. To obtain information
about $M$, one can consider the base change homomorphisms
$$
\Delta_I(M) = \Z[G^{n}] \otimes_{\Z[G]} M
$$
where the algebra homomorphism $\Delta_I:\Z[G] \to \Z[G^{n}]$ is induced 
by the homomorphism
 $I:G \to G^{n}$.  The statement that certain tensor products
over $\Z[G^{n}]$ of these $\Delta_I(M)$ have trivializations puts a constraint
on $M$, which in some cases one might hope would force $M$ to have a trivial
class in $\Cl(\Z[G])$.  

The particular tensor product which we consider is
the  so-called $\Theta_n$-product defined by
\begin{equation*}
\label{thetanowdef}
\Theta_n(M) = \bigotimes_{I \subset \{1,\ldots, n\}} \Delta_I(M)^{(-1)^{n-\#I}}.
\end{equation*}
We define an ``$n$-cubic structure" on $M$ to be a 
trivialization of $\Theta_n(M)$ that satisfies certain ``cubic" compatibility  conditions.
(The terminology and notation come from the theorem of the cube for line bundles
on abelian varieties and 
the associated theta functions; this classical case corresponds to $n=3$.)

We will like to apply the above to $M=\det \R\Gamma(X,\O_X)^{\otimes 2}$ and $n=d+2$.
To  extract the maximal amount of
information about $M$, it is important to consider in addition to $M$ a generator $s$ 
of $M_\Q = \Q \otimes_\Z M$ as a $\Q[G]$-module.  
Such an $s$ gives rise to generators
$\Delta_I(s) = 1 \otimes_{\Z[G]} s$ and $\Theta_{n}(s)=\prod_{I}\Delta_I(s)^{(-1)^{n-\#I}}$ of 
$\Delta_I(M_\Q)$ and $\Theta_{n}(M_\Q)$ respectively.

In general, this $M$ does not support an $n$-cubic structure. However, $M_\Q$ does 
because the cover $X_\Q\to Y_\Q$ is unramified.    When 
$d= {\rm dim}(X_\Q) = 1$ this follows from the Deligne-Riemann-Roch theorem [De]; for
general $d$ it follows from work of Breen [Br], Ducrot [Du] and Pappas [P3]. 
The resulting trivialization of $\Theta_{n}(M_\Q)$ enables us to
regard  $\Theta_{n}(s)$ as an element $c_n(s)\in \Q[G^n]^*$. One 
can view $c_n(s)$ as a kind of secondary class, analogous to 
the Bott-Chern secondary characteristic class arising from discrepancies 
between metrics on vector bundles lying in an exact
sequence of vector  bundles. 

Now suppose that we choose a generator
$s_v$ of $\Z_v\otimes_\Z M$ as a $\Z_v[G]$-module at each finite 
place $v$ of $\Q$. The above procedure (performed now over $\Q_v$) allows us to obtain
elements $c_n(s_v)\in \Q_v[G^n]^*$. The elements $c_n(s_v)$
are our substitutes of the local (abelian) resolvents that appear 
in Fr\"ohlich's theory of resolvents for actions of finite groups on rings
of algebraic integers.  Our main calculation
amounts to showing that, under certain assumptions, we can compute
$c_n(s_v)$ up to an element of $\Z_v[G^n]^*$ using resolvent theory on the codimension $1$ points of $Y$ and a 
localized Riemann-Roch theorem. (Note that in the main text we mainly work with the $\Z[G^n]$-ideal
in $\Q[G^n]$ given by the id\`ele $(c_n(s_v))_v$; it is denoted by $E$
in \S 3.)

Our results are actually expressed in terms of a functor
\begin{equation*}
\Theta_n: {\cal Pic}^\eta(\Z[G])\to {\cal C}_\Z(G;n); \quad (N, t)\mapsto (\Theta_n(N), \Theta_n(t))
\end{equation*}
where ${\cal Pic}^\eta(\Z[G])$ is the Picard category (see [SGA4, XVIII]) of locally free rank one $\Z[G]$-modules $N$
with a generator $t$ of $N_\Q$ and ${\cal C}_\Z(G;n)$ is a similar Picard category of locally free rank one $\Z[G^n]$-modules
where the morphisms are restricted according to certain  ``cubic" conditions
(which will be specified in detail in the next section).  
Considering isomorphism classes of objects in Picard categories leads
 to a homomorphism
\begin{equation*}
\Theta_n:\Pic(\Z[G]) \to C_\Z(G;n)
\end{equation*}
whose kernel  is the group  of classes in $\Pic(\Z[G])$ having an $n$-cubic structure.  
By the results in [P2] this kernel is small (see Theorem \ref{2.7}; it is actually trivial when $n\leq 5$). The main result of this paper, Theorem \ref{mainthm}, 
gives an explicit  ``branch locus" formula for $\Theta_n(M)$ when
$n=d+2$, $M=\det \R\Gamma(X,\O_X)^{\otimes 2}$ and the ramification of $\pi:X \to Y = X/G$
is domestic, in the sense that it is supported away from
the prime divisors of $\# G$.

Let us conclude with some remarks.
We first note that the approach outlined above allows us to
provide a new explanation 
of the ``cubic"  conditions on the trivializations of $\Theta_n(M)$
which were first considered by Breen ([Br]) when $n=3$. We show 
that Breen's conditions are among those satisfied by all ``ratios"
$\Theta_n(s)/\Theta_n(s')\in \Q[G^n]^*$ where $s$ and $s'$ are two generators
of $\Q\otimes_\Z M$ (see Remark  \ref{rem1} (b)). We can consider a bigger set of conditions and obtain variants 
of the notion of cubic structure
and of the functor $\Theta_n$ above (see Remarks \ref{rem1} (d) and \ref{remvar}). In fact, by taking 
the maximal set of such conditions
we can obtain a possibly more natural variant of the notion of cubic structure.
When $n=2$ the homomorphism $\Theta_2$ is very closely related to the 
functor ${\rm rag}$ of McCulloh [McC]. Our formulas then essentially specialize 
to the formulas in loc. cit.  

A major remaining problem is to understand the nature
of the ``cubic relations" when $G$ is non abelian. Some progress has been made
on that problem by the second author when $n=3$ ([P5]). (This pertains 
to covers of arithmetic surfaces.) Also, since the paper 
was submitted there has been some progress by the authors in developing 
a refined Riemann-Roch theorem for  
non-abelian actions on arithmetic surfaces. However, we have not, as yet, been able to establish a comparably precise result about Euler characteristics even in the surface case. The situation for $n> 3$ remains
completely unclear at the moment.  One of the main themes of [P4-5] is that, 
because the $\gamma$-filtration 
on the Grothendieck group of vector bundles on $Y=X/G$ terminates after the $d+1$-th step,
there are relations between the Chern classes of vector bundles
over $Y$ which are obtained from various representations
of $G$ using the cover $X\to Y$.   Such relations alone can sometimes be used to deduce 
information about an equivariant Euler characteristic  $\overline{\chi}^P(X, \F)$.
In this paper we have taken advantage of the fact that the $\gamma$-filtration on the representation ring of an abelian group 
has a very simple form.  The $\gamma$-filtration on the representation ring
of a general finite group has no simple and uniform explicit description;
this is at least one reason why it is hard to 
extend our results to arbitrary finite  groups. 

\smallskip

\noindent {\bf Acknowledgement:}{ We would like to thank L. Illusie, N. Katz, L. Merel and P. Sarnak for useful conversations and the referee for helpful comments.}

\section{Cubic structures and categories}\label{cube} 
\subsection{} 

Let $G$ be a finite
abelian group, and let $R$ be a Dedekind ring with fraction field $K$.
Let $n\geq 2$ be an integer. 
We will denote by $G^D_R=\Spec(R[G])$ the Cartier dual of the constant group scheme given by $G$.
We start by rephrasing the notion of $n$-cubic structure on a line bundle $\L=\widetilde M$
over $G^D_R$ (see [P2, \S 3]) in terms of the corresponding locally free rank one $R[G]$-module $M=\Gamma(G^D_R, \L)$.

Suppose $r$ and $s$ are non-negative  integers.  Denote by $G^r $ the direct sum of
$r$ copies of $G$, where $G^{0}$ is taken to mean the group with one element. For 
$\mathcal{I} \in C_{r,s} = \mathrm{Hom}(G^{r},G^{s})$, we will  denote by 
$\Delta_{\mathcal{I}}:R[G^r] \to R[G^s]$ the induced $R$-algebra
homomorphism. These homomorphisms have the property that
\begin{equation}
\label{eq:complaw}
\Delta_{\mathcal{I'} \circ \mathcal{I}} = \Delta_{\mathcal{I'}} \circ \Delta_{\mathcal{I}}
\end{equation}
if $\mathcal{I'} \in C_{s,s'}$ for some integer $s'$. Also, we will denote by 
$\Delta^D_{\mathcal I}: (G^D_R)^s\to (G^D_R)^r$ the $R$-scheme morphism
given by $\Delta_{\mathcal I}$. 

Let $\Sigma_{r,s} = {\bf Z}[C_{r,s}]$. Define the bilinear map
\begin{equation}
\label{eq:Smult}
\Sigma_{r,s} \times \Sigma_{r',s'} \to \Sigma_{r',s}
\end{equation}
to be  induced by the composition of homomorphisms $C_{r,s} \times C_{r',s'}\to C_{r',s}$ 
when $s' = r$ and to be the zero map otherwise.
These maps make  $\Sigma = \bigoplus_{r,s \ge 0} \Sigma_{r,s}$ into a ring without unit.

Now let $M$ be a locally free rank one $R[G^r]$-module. We will denote 
by $M^{-1}$ the dual locally free rank one $R[G^r]$-module ${\rm Hom}_{R[G^r]}(M, R[G^r])$.
Let $t$ be a generator of $M_K = K\otimes_RM$ as a free rank one
$K[G^r]$ -module.   
Define 
\begin{equation}\label{Deldefine}  \Delta_{\mathcal{I}}\left(  M\right)  = 
R[G^s]\otimes_{\Delta_{\mathcal{I}}, R[G^r]} M  \ ,\qquad    \Delta_{\mathcal{I}}(t) = 1 \otimes t
\in  
 \Delta_{\mathcal{I}}\left( 
M\right)_K
\end{equation}
so that  $ \Delta_{\mathcal{I}}\left(  M\right)$ is a locally free rank one $R[G^s]$-module, and
$\Delta_{\mathcal{I}}(t)$ is a generator of $ \Delta_{\mathcal{I}}\left(  M\right)_K$.
In general, suppose $z = \sum_{\mathcal{I} \in C_{r,s}}
z({\mathcal{I}})
\cdot \mathcal{I} \in
\Sigma_{r,s}$ for some integers
$z(\mathcal{I})$.  We define
\begin{equation}
\label{eq:zdef}
\Delta_z(M) = \bigotimes_{\mathcal{I}} \Delta_{\mathcal{I}}(M)^{ z(\mathcal{I})}
\quad \mathrm{and}\quad
\Delta_z(t) = \bigotimes_{\mathcal{I}} \Delta_{\mathcal{I}}(t)^{z(\mathcal{I})}\ 
\end{equation}
where the tensor products are over $R[G^s]$, resp. $K[G^s]$.  We also have
a homomorphism of multiplicative groups
$\lambda_z:K[G^r]^*\rightarrow K[G^s]^*$
defined by
\begin{equation}
\label{eq:zform}
\lambda_z(\alpha) = \prod_{\mathcal{I}} \Delta_{\mathcal{I}}(\alpha)^{z(\mathcal{I})}.
\end{equation}

If $z' \in \Sigma_{s,s'}$ for some integer $s'$, then 
the composition law (\ref{eq:complaw}) gives a canonical isomorphism
\begin{equation}
\label{eq:nicecomp}
\Delta_{z'} \circ \Delta_z \simeq \Delta_{z'\cdot z}\ .
\end{equation}
In particular, if $z' \cdot z = 0$ in the ring $\Sigma$, then there is a canonical
isomorphism
\begin{equation}
\label{eq:trivial}
\Delta_{z'} (\Delta_z(M)) \simeq R[G^{s'}].
\end{equation}
The corresponding isomorphism over $K$ sends $\Delta_{z'}(\Delta_z(t))$ to $1 \in K[G^{s'}]$.

A subset $I \subset \{1,\ldots , n\}$ determines
a homomorphism 
$I \in C_{1,n} = \mathrm{Hom}(G,G^{n})$ 
by defining $I(g)$ for $g \in G$ to have $i^{th}$ component $g$
for $i \in I$ and $i^{th}$ component the identity element $e$ of $G$ if $i \notin I$.  In this way, we will
view subsets $I$ of $\{1,\ldots , n\}$ as elements of $C_{1,n}$.
Notice that the scheme morphism $m_I=\Delta_I^D: (G^D_R)^n=\Spec(R[G^n])\to G^D_R=\Spec(R[G])$
induced by $\Delta_I$ is given on points by $m_I(x_1,\ldots, x_n)=\sum_{i\in I}x_i$.
(Note that we use additive notation for the group operation.)
Now consider the element of $\Sigma_{1,n}$ given by
\begin{equation*}
s_n=\sum_{I\subset \{1,\ldots, n\}}(-1)^{n-\#I}\cdot I\ .
\end{equation*}
If $M$ is a locally free rank $1$ $R[G]$-module, for simplicity, 
we will denote $\Delta_{s_n}(M)$ and $\Delta_{s_n}(t)$
by $\Theta_n(M)$ and $\Theta_n(t)$ respectively. Hence, we have
\begin{equation*}
\Theta_n(M)=\bigotimes_{I \subset \{1,\ldots , n\}}\Delta_{I}(M)^{ (-1)^{n-\#I}}. 
\end{equation*}
\smallskip

We now define three conditions on an element $a \in K[G^n]^*$. 
\smallskip

1.  Consider the trivial
homomorphism $e: g \to e$ in  $C_{n,0}={\rm Hom}(G^n, G^0)\subset \Sigma_{n,0}$.  
Since $n \ge 1$,  $e \cdot
s_n = 0$ in $\Sigma_{1,0}$.   The homomorphism $\lambda_e: K[G^n]^* \to K^*$ is the one induced by the 
trivial character of $G^n$. An $a \in K[G^n]^*$ such that
$\lambda_{e}(a) = 1$ will be said to be {\sl rigid}.  
\medskip

2. Suppose  that $\sigma$ is a permutation of $\{1,\ldots, n\}$.  We then have
an element
$$\mathcal{I}_\sigma :(g_1,\ldots, g_n) \to (g_{\sigma(1)},\ldots, g_{\sigma(n)})$$ of $C_{n,n}$.
Let ${\rm id}$ be the identity permutation of $\{1,\ldots, n\}$, and let
$z_\sigma = \mathcal{I}_\sigma - \mathcal{I}_{\rm id}$ in $\Sigma_{n, n}$.  One checks readily that
$z_\sigma \cdot s_n = 0$ in $\Sigma_{1, n}$.  An element $a \in K[G^n]^*$ such that
$\lambda_{z_\sigma}(a) = 1$ for all permutations $\sigma $ of $\{1,\ldots, n\}$ will 
be said to be {\sl symmetric}.
\medskip

3. Define four elements
of $C_{n, n+1}=\mathrm{Hom}(G^n,G^{n+1})$ by
\begin{eqnarray*}
\label{eq:assocrels}
 \mathcal{I}_0 &:& (g_1,\ldots, g_n) \to (g_1,g_1,g_2\ldots, g_n)  \cr
\mathcal{I}_1 &:& (g_1,\ldots, g_n) \to (e,g_1,g_2\ldots,g_n)  \cr
 \mathcal{I}_2 &:& (g_1,\ldots, g_n) \to (g_1,g_2,e,\ldots, g_n) \cr
\mathcal{I}_3 & :& (g_1,\ldots, g_n) \to (g_1,g_2,g_2,\ldots, g_n) 
\end{eqnarray*}
Set
$z = \mathcal{I}_0 - \mathcal{I}_1 + \mathcal{I}_2 -\mathcal{I}_3$ in $\Sigma_{n, n+1}$.
One has $z \cdot s_n = 0$ in $\Sigma_{1,n+1}$ since $n \ge 2$.  Elements $a \in K[G^n]^*$
such that $\lambda_z(a) = 1$ will be said to satisfy the {\sl cocycle condition}.  

\begin{Definition}
{\rm ([P2, 3.3]) An element $a\in K[G^n]^*$ which is rigid, symmetric and satisfies the cocycle condition
(see 1,2 and 3 above) will be said to be {\sl $n$-cubic.} } 
\end{Definition}

\begin{Definition}\label{cubicdef}
{\rm Suppose that $R'$ is a subring of $K$ that contains $R$.
An $n$-cubic structure on a locally free rank $1$ $R'[G]$-module $M$
is a trivialization 
$$
c: R'[G^n]\xrightarrow{\sim} \Theta_n(M)
$$
such that for any choice of generator $t$ of $M\otimes_{R'}K$ as a $K[G]$-module,
we have $c(1)=a\cdot \Theta_n(t)$ with $a\in K[G^n]^*$ which is $n$-cubic.}
\end{Definition}

\begin{Remarknumb}\label{rem1}{\rm a) In the above definition, we are mainly interested 
in the cases $R'=R$ or $R'=K$.

b) Suppose that $t'$ is another choice of $K[G]$-generator of $M\otimes_{R'}K$
so that $t'=\alpha\cdot t$ with $\alpha\in K[G]^*$. Then 
$\Theta_n(t')/\Theta_n(t)=\lambda_{s_n}(\alpha)$. Now if $z$ is as in 
conditions $1$, $2$ or $3$ above, we have (as in (\ref{eq:trivial}))
$\lambda_{z}(\lambda_{s_n}(\alpha))=\lambda_{z\cdot s_n}(\alpha)=\lambda_0(\alpha)=1$.
Therefore, the element $\lambda_{s_n}(\alpha)$ is $n$-cubic.
Hence, it is enough to check the property of the definition, for a single
choice of $K[G]$-generator of $K\otimes_{R'}M$.

c) We can readily see that the notion of $n$-cubic element coincides with the
corresponding notion as defined in [P2]. There the conditions
are expressed in terms of characters of $G$,
i.e of points of $G^D_R$. Also, $n$-cubic structures
on $M$, as defined above, uniquely correspond to $n$-cubic structures
(as defined in [P2, \S 3.a])
on the corresponding line bundle $\widetilde M$ over $G^D_{R'}$. 

d) We can consider the following variants of the notion of $n$-cubic structure:
Let $V$ be a set of elements $z\in \bigcup_{t\geq 0}\Sigma_{n,t}$ that satisfy
$z\cdot s_n=0$. Then, we can consider ``$V$-cubic elements": by definition, these are elements
$a\in K[G^n]^*$ such that $\lambda_z(a)=1$ for all $z\in V$. The notion of ``$V$-cubic structure"
is defined as in Definition \ref{cubicdef} above by replacing ``$n$-cubic" by ``$V$-cubic". 
Notice that the argument in Remark (b) shows that for $a\in K[G]^*$, 
$\lambda_{s_n}(a)$ is $V$-cubic. Therefore, Remark (b) applies also
to $V$-cubic structures.}
\end{Remarknumb}

\subsection{ }\label{cubecat}
The reader is referred to [SGA4, XVIII] for the definition
of a Picard category. We now define two such categories:

The category ${\cal Pic}^\eta_R(G)$ with objects pairs $(M, t)$ 
where $M$ is a locally free rank $1$ $R[G]$-module and $t$ a
generator of $K\otimes_RM$ as a $K[G]$-module, and morphisms
$\phi: (M, t)\to (M', t')$ given by $R[G]$-isomorphisms $\phi: M\to M'$
(with no condition on the generators $t$ and $t'$).

The category ${\cal C}_R(G;n)$ with objects pairs $(P, \gamma)$ 
where $P$ is a locally free rank $1$ $R[G^n]$-module and $\gamma$ a
generator of $K\otimes_RP$ as a $K[G^n]$-module, and morphisms
$\psi: (P, \gamma)\to (P', \gamma')$ given by $R[G^n]$-isomorphisms 
$\psi: P\to P'$ such that $({\rm id}_K\otimes_R\psi)(\gamma)=a\cdot \gamma'$ where $a$
is an $n$-cubic element of $K[G^n]^*$.

We can see that both ${\cal Pic}^\eta_R(G)$ and ${\cal C}_R(G; n)$ are strictly commutative Picard categories
with product defined by $(M, t)\otimes (M', t')=(M\otimes_{R[G]}M', t\otimes t')$
and $(P, \gamma)\otimes (P', \gamma')=(P\otimes_{R[G^n]}P', \gamma\otimes \gamma')$.
The group of isomorphism classes of ${\cal Pic}^\eta_R(G)$
is the Picard group ${\rm Pic}(R[G])$; we will denote by $C_R(G;n)$ the (abelian)
group of isomorphism classes of objects of ${\cal C}_R(G;n)$.

\begin{lemma}
There is an additive functor $\Theta_n: {\cal Pic}^\eta_R(G)\to {\cal C}_R(G;n)$
given by 
\begin{equation*}
\Theta_n((M, t))=(\Theta_n(M), \Theta_n(t)).
\end{equation*}
\end{lemma}

\begin{Proof} It follows from the argument in Remark \ref{rem1} (b) that
the functor $\Theta_n$ is well-defined. The rest is left to the reader.\endproof
\end{Proof}

The functor $\Theta_n$ induces a group homomorphism
\begin{equation}\label{eq:theta}
\Theta_n: \Pic(R[G])\to C_R(G; n).
\end{equation}
When the integer $n$ is fixed from the context 
we will simply write $\Theta$ instead of $\Theta_n$.

\begin{Remarknumb} \label{remkernel}
{\rm  Notice that by the definitions, isomorphisms $(R[G^n], 1)\xrightarrow{\sim} \Theta((M, t))$
in the category ${\cal C}_R(G;n)$ correspond to $n$-cubic structures on the $R[G]$-module $M$.
Hence, the kernel of (\ref{eq:theta}) is the set (which is actually a group under tensor product) of isomorphism
classes of locally free rank $1$ $R[G]$-modules which support an $n$-cubic structure. 
 }
\end{Remarknumb}
\begin{Remarknumb}\label{remvar}
{\rm Suppose $V$ is a set as in Remark \ref{rem1} (d). Then there are obvious variants 
${\cal C}_R(G;V)$ and $\Theta_V:   {\cal Pic}^\eta_R(G)\to {C}_R(G;V)$,
of ${\cal C}_R(G;n)$ and  $\Theta_n$ respectively. They are given by replacing 
``$n$-cubic element" by ``$V$-cubic element" in the above definitions. 
}
\end{Remarknumb}

\subsection{}\label{ideles}
 Assume that $R$ is the ring of integers of a number field $K$.
For $m\geq 0$, let us denote by ${\bf A}^*_{f, K[G^m]}$ the finite id\`eles of $K[G^m]$. Then
\begin{equation*}
{\bf A}^*_{f, K[G^m]}=\prod_{v}{}^{'}K_v[G^m]^*
\end{equation*}
where the (restricted) product is over all 
finite places $v$ of $K$. A finite id\`ele $(a_v)_v$ for $K[G^m]$ gives a locally free
rank $1$ $R[G^m]$-module $Q((a_v)_v)$ by
\begin{equation*}
Q((a_v)_v)= \cap_{v}(  R_v[G^m]a_v\cap K[G^m])\subset K[G^m].
\end{equation*}
The ``generic fiber" $Q((a_v)_v)\otimes_RK=K[G^m]$ has $1$ as a distinguished generator
and we have $Q((a_v)_v)\otimes_{R}R_v= R_v[G^m]a_v$.
Therefore, the id\`ele $(a_v)_v$ provides us with an object $(Q((a_v)_v), 1)$  
of the category ${\cal C}_R(G;m)$ if $m\geq 2$ or of ${\cal Pic}^\eta_R(G)$ 
if $m=1$. 

For each finite place $v$ of $K$ fix an algebraic closure $\bar K_v$ of the completion $K_v$.
Let us denote by ${\rm Ch}(G^n)_v$ the (additive) ring of virtual characters of $G^n$ 
with values in $\bar K_v$. There is a group isomorphism
\begin{equation}
K_v[G^{n}]^*\xrightarrow{\sim} {\rm Hom}_{{\rm Gal}(\bar K_v/K_v)}({\rm Ch}(G^n)_v, \bar K_v^*)
\end{equation} 
given by $a_v\mapsto (\phi\mapsto \phi(a_v))$ (cf. [Fr, II \S 1]). 

Notice that for every element $\mathcal I$ of $C_{r,s}$ 
we can obtain an additive map
$\Delta_{\mathcal I}^D: {\rm Ch}(G^s)_v\to {\rm Ch}(G^r)_v$ 
by composing a character of $G^s$ with the ring homomorphism
$\bar K_v[G^r]\to \bar K_v[G^s]$ given by $\Delta_{\cal I}$ and then extending   linearly
to virtual sums of such characters. This definition generalizes to elements 
$z=\sum_{{\mathcal I}\in C_{r,s}}z({\cal I})\cdot {\cal I}$
of $\Sigma_{r,s}$ by linearity; we obtain
additive maps $\Delta_z^D: {\rm Ch}(G^s)_v\to {\rm Ch}(G^r)_v$.

There is a commutative diagram
\begin{equation}\label{idel1}
\begin{CD}
\prod'{\rm Hom}_{{\rm Gal}(\bar K_v/K_v)}({\rm Ch}(G)_v, \bar K_v^*)@>\sim>> {\bf A}^*_{f, K[G]}@>Q>> \Pic(R[G])\\
@V \prod'(\Theta_n^D)^* VV  @V \prod'\lambda_{s_n} VV @VV\Theta_n V\\
\prod'{\rm Hom}_{{\rm Gal}(\bar K_v/K_v)}({\rm Ch}(G^n)_v, \bar K_v^*)@>\sim>> {\bf A}^*_{f, K[G^n]}@>Q>> C_R(G;n)
\end{CD}
\end{equation}
where $(\Theta_n^D)^*$ is the map dual to $\Theta^D_n:=\Delta^D_{s_n}: {\rm Ch}(G^n)_v\to {\rm Ch}(G)_v$.
Explicitly, we have
\begin{equation}\label{naugment}
\Theta^D_n(\phi_1\otimes\cdots \otimes\phi_n)=
\sum_{I\subset\{1,\ldots, n\}}(-1)^{n-\#I}\prod_{i\in I}\phi_i=\prod_{i=1}^n(\phi_i-1).
\end{equation}
Notice that the image of $\Theta_n^D$ is the $n^{th}$ power of the augmentation ideal 
of ${\rm Ch}(G)_v$.

\subsection{}\label{func}
Suppose that $\phi: G\to H$ is a homomorphism of finite abelian groups; this induces
ring homomorphisms $ R[G]\to R[H]$, $\phi^{(n)}: R[G^n]\to R[H^n]$. If $c\in K[G^n]^*$
is an $n$-cubic element then so is its image $\phi^{(n)}(c)\in K[H^n]^*$. Tensoring with $R[H]$ over $R[G]$, resp.  
with $R[H^n]$ over $R[G^n]$, gives additive functors $ {\cal Pic}^\eta_R(G)\to {\cal Pic}^\eta_R(H)$,
resp. $ {\cal C}_R(G; n)\to {\cal C}_R(H; n)$. The corresponding diagram
\begin{equation}
\begin{CD}\label{functoriality}
 {\rm Pic}(R[G])@>\Theta>> C_R(G;n)\\
@VVV @VVV\\
 {\rm Pic}(R[H])@>\Theta>> C_R(H;n)
\end{CD}
\end{equation}
commutes. There is also a commutative diagram
\begin{equation}
\begin{CD}\label{idel2}
\prod'{\rm Hom}_{{\rm Gal}(\bar K_v/K_v)}({\rm Ch}(G^n)_v, \bar K_v^*)@>\sim>> {\bf A}^*_{f, K[G^n]}@>Q>>
C_R(G;n)\\ @V \prod'{\rm C}^*_\phi VV  @V \prod'\phi^{(n)} VV @VV\phi  V\\
\prod'{\rm Hom}_{{\rm Gal}(\bar K_v/K_v)}({\rm Ch}(H^n)_v, \bar K_v^*)@>\sim>> {\bf A}^*_{f, K[H^n]}@>Q>> C_R(H;n)
\end{CD}
\end{equation}
where ${\rm C}_\phi^*$ is  dual to the map ${\rm C}_\phi: {\rm Ch}(H^n)_v\to {\rm Ch}(G^n)_v$
given by composing  characters with $\phi^n: G^n\to H^n$.

\subsection{}
Let us now assume that $R=\Z$. Let ${\mathcal M}_G$ be the normalization of $\Z[G]$ in $\Q[G]$.
Then ${\mathcal M}_G$ is the maximal  order in $\Q[G]$. Tensoring with ${\mathcal M}_G$ over $\Z[G]$
induces a homomorphism $\Pic(\Z[G])\to \Pic({\mathcal M}_G)$. 
Denote its kernel by ${\rm D}(\Z[G])$; by Rim's theorem it is trivial when $G$ is of prime order.

Now recall that the $k$-th Bernoulli number $B_k$ is defined by the power series identity:
 $ {{t}/{(e^t-1)}=\sum^{\infty}_{k=0}B_k{t^k}/{k!}}$. For a prime $p$ we denote by $|\ |_p$
the usual $p$-adic absolute value with $\displaystyle{|p|_p=p^{-1}}$.
Let $e(1)=1$ and for $k\geq 2$ let us set 
\begin{equation*}
e(k)=\begin{cases}
\ \ {\rm numerator\,}({B_{k}}/{k})&,\ \text{if $k$ is even},\\
\displaystyle{\prod_{p,p|h^+_p }}|\#\Kr_{2k-2}(\Z)|_p^{-1}&,\ \text{if $k$ is odd},
\end{cases}
\end{equation*}
where $\Kr_{2k-2}(\Z)$ is the Quillen $\Kr$-group 
(a finite group for $k>1$) and we have $h^+_p=\# {\rm Cl}(\Z[\zeta_p+\zeta_p^{-1}])$. 
Note that, according to Vandiver's conjecture, $p$ does not divide $h^+_p$,
which implies $e(k)=1$ for $k$ odd.
The following is essentially one of the main results of [P2]
(or [P1] when $n=3$).

\begin{thm}\label{2.7}  Suppose that $R=\Z$ and that $\Theta_n: 
{\Pic}(\Z[G]) \to  {C}_{\Z}(G; n)$
is the homomorphism (\ref{eq:theta}).
\begin{enumerate}
\item[a.] ${\rm ker}(\Theta_n)$ is annihilated by $\prod_{k=1}^{n-1}\prod_{p|e(k)}|\#G|_p^{-1}$.

\item[b.] If $2\leq n\leq 5$, ${\rm ker}(\Theta_n)$ is trivial.

\item[c.] If all the prime divisors of $\#G$ are greater than or equal to $n$
and satisfy Vandiver's conjecture then ${\rm ker}(\Theta_n)\subset {\rm D}(\Z[G])$.
In particular, if $\#G=p$ is a prime number $\geq n$ which satisfies Vandiver's conjecture 
then ${\rm ker}(\Theta_n)$ is trivial. 
\end{enumerate}
\end{thm}
\begin{Proof}   Recall that by Remark \ref{remkernel} elements in the kernel 
${\rm ker}(\Theta_n)\subset {\rm Pic}(\Z[G])$
are represented by invertible sheaves that support an $n$-cubic structure. Therefore,
 parts (a) and (b) are given by  [P2, Theorem 1.1]. 
Part (c)  follows
from [P2, Theorem 1.2].\endproof 
\end{Proof}

\section{The main calculation and result} \label{mainsec}
We continue to assume that $R$ is a Dedekind ring with
field of fractions $K$. Suppose that $h: Y\to \Spec(R)$ is a regular flat projective scheme,
equidimensional of (absolute) dimension $d+1$. Let $\pi: X\to Y$ be a $G$-cover
where $G$ is a finite abelian group. By definition, this means that $X$ supports a (right) action of
$G$,  $\pi$ is finite and identifies $Y$ with the quotient $X/G$,
and that $\pi$ is generically on $Y$ a $G$-torsor. Denote by $f:X\to \Spec(R)$ the structure morphism.
We assume  that $X$ is normal,  in addition to the following hypothesis:  
\smallskip
 
(T) The action of $G$ on $X$ is {\sl tame}, i.e
for every point $x$ of $X$, the order of the inertia 
subgroup $I_x$ of $x$ is relatively prime to the
residue field characteristic of $x$.
In addition, we assume that the cover $\pi_K: X_K\to Y_K$ over $\Spec(K)$ 
is unramified
(i.e a $G$-torsor; then $X_K$ is also regular).
This last condition follows from the assumption on the tameness of the action, if
$\Spec(R)\to \Spec(\Z)$ is surjective. 
\smallskip

\begin{Remarknumb} \label{abhya}
  {\rm a)  It follows from the above assumptions that $\pi: X\to Y$ is flat.
To show this, let $x$ and $y$ be points of $X$ and $Y$ 
with $\pi(x)=y$ and let  $\hat Y$, resp. $\hat X$, be the spectrum of the strict 
henselization of the local ring
$\O_{Y, y}$, resp. $\O_{X,x}$. By descent, it is enough to show
that $\hat \pi: \hat X\to \hat Y$ is flat. By [Ra] we have 
\begin{equation}\label{inducedI}
X\times_Y \hat Y\simeq \hat X\times^{I_x} G:=(\hat X\times G)/I_x
\end{equation}
as schemes with $G$-action. Therefore, $\hat X/I_x\simeq \hat Y$.
By our assumptions, $\#I_x$ is relatively prime to the residue characteristic,
$\hat Y$ is regular and $\hat X$ is normal. The argument in [Ro1] now
shows that $\hat \pi$ is flat.

b) To construct examples in which hypothesis (T) holds, one can start with 
a regular projective model $Y$ of $Y_K$ and a $G$-torsor $\pi_K:X_K \to Y_K$.
One then considers whether the normalization $X$ of $Y$ in the total quotient
ring $K(X_K)$ of $X_K$
is tame over $Y$, using for example Abhyankar's Lemma [GM, 2.3.2].  Suppose,
for instance, that  $Y'$ is a regular curve over $\Spec(R[1/N])$ for some integer $N$
with $(N, \#G)=1$, and that $\pi':X' \to Y'$ is a $G$-torsor.  Resolution of 
a $\Spec(R)$-model for $Y'$ 
then leads to a $G$-cover $X\to Y$ over $\Spec(R)$ which satisfies
all of our assumptions and has general fiber $\pi'_K$. In fact, this cover has at worst ``domestic ramification"
(see \S \ref{mainthmpar}).} 
\end{Remarknumb}
\smallskip

Under the above conditions, the sheaf $\pi_*(\O_X)$ of $\O_Y[G]$-modules
on $Y$ is locally $\O_Y[G]$-free of rank $1$. (Since by Remark \ref{abhya} (a) $\pi$ is flat,
this follows from [CEPT1, Theorem 2.6 and Proposition 2.7].)
Hence, we may think of $\pi_*(\O_X)$ as an invertible sheaf  (line bundle)
on the scheme $Y\times_{\Spec(R)}G^D_R=G^D_Y$.  For simplicity, we will denote this invertible sheaf by
$\L$. Similarly, if $\Gg$ is a locally free coherent $\O_Y$-sheaf on $Y$ we can consider the $\O_X$-sheaf 
$\F=\pi^*(\Gg)$ with compatible $G$-action and regard $\pi_*(\F)$ as a locally free coherent 
$\O_{G^D_Y}$-sheaf on $G^D_Y$.
We then have $\pi_*(\F)\simeq\Gg\otimes_{\O_Y}\pi_*(\O_X)=\Gg\otimes_{\O_Y}\L$. 
In general, if $\F$ is any  locally free coherent $G$-sheaf on $X$
(i.e a locally free coherent $\O_X$-module with $G$ action compatible
with the action of $G$ on $X$) then $\pi_*(\F)$ 
can also be thought of as a locally free coherent $\O_{G^D_Y}$-sheaf on $G^D_Y$. 

Denote by $\ti h: G^D_Y\to G^D_R$ the base change of $h: Y\to \Spec(R)$ by $G^D_R\to \Spec(R)$.
Suppose that $\Hh$ is a locally free coherent $\O_{G^D_Y}$-module.
The total derived image ${\R}\ti h_*(\Hh)$ 
in the derived category of complexes of sheaves of $\O_{G^D_R}$-modules on $G^D_R$ which are bounded below
is ``perfect" (i.e  it is locally quasi-isomorphic to a 
bounded complex of   locally free coherent $\O_{G^D_R}$-modules, see [SGA6, \S III]).  Hence, by 
[KMu], we can associate to ${\R}\ti h_*(\Hh)$ a graded invertible 
sheaf  \begin{equation*}
{\det}_* {\R}\ti h_*(\Hh) =(\det {\R}\ti h_*(\Hh), {\rm rank}({\R}\ti h_*(\Hh))\, {\rm mod}\, 2)
\end{equation*}
on $G^D_R$ (the ``determinant of cohomology"). In what follows, we will mostly consider
situations in which the second term   ${\rm rank}({\R}\ti h_*(\Hh))\, {\rm mod}\, 2$
(a Zariski locally constant $\Z/2\Z$-function on $G^D_R$)
is trivial. Then we will call $\det {\R}\ti h_*(\Hh)$ ``the determinant 
of cohomology".

Suppose now that $R$ is the ring of integers
of a number field. The tameness assumption allows us  to define 
the equivariant projective Euler characteristic $\chi^P(X, \F)$ in the Grothendieck group
of finitely generated projective $R[G]$-modules $\Kr_0(R[G])$
([CEPT1, Theorem 8.3], see also [C]).   We denote 
by $\bar\chi^P(X,\F)$ the class of $\chi^P(X,\F)$ in the quotient 
$\Cl(R[G]):=\Kr_0(R[G])^{\rm red}=\Kr_0(R[G])/\pm\{\hbox{\rm free $R[G]$-modules}\}$ 
(the ``projective class group of $R[G]$"). Since $G$ is abelian
and $R$ has Krull dimension $1$, there is a natural identification
\begin{equation*}
\Pic(G^D_R)=\Pic(R[G])=\Cl(R[G]).
\end{equation*}
Under this identification we have ([P1, \S 2.c])
\begin{equation}\label{identity}
[\det {\R}\ti h_*(\pi_*(\F))]=\bar\chi^P(X,\F).
\end{equation}

\subsection{Resolvent theory over $Y$} 
Let 
$\chi: G\to \Gamma(S',\O^*_{S'})=\Gm(S')$
be a character with $S'\to \Spec(R)$ an $R$-scheme; then
$\chi$ corresponds to an $S'$-point of $G^D_R$.
Set $\pi': X'=X\times_R S'\to Y'=Y\times_R S'$. Under the assumption (T),
$\pi'$ identifies $Y'$ with the quotient $X'/G$ ([CEPT1, Theorem 5.1], or [KM, Proposition A7.1.3]
combined with (\ref{inducedI})).
Recall $G$ acts on $(\pi')_*(\O_{X'})=\pi_*(\O_X)\otimes_R\O_{S'}$ via its action on $\O_X$.
Consider the subsheaf of $(\pi')_*(\O_{X'})$
given by local sections $b \in \O_{X'}=\O_X\otimes_R\O_{S'}$ that satisfy
$g\cdot b = \chi(g)b$.
We will denote this subsheaf by $\O_{X',\chi}$. We can see that, under our assumptions, $\O_{X',\chi}$
is an invertible $\O_{Y'}$-sheaf on $Y'$. (In accordance with the classical theory of Fr\"ohlich
one might call $\O_{X',\chi}$ a {\sl resolvent} invertible sheaf over $Y'$.)
 
\begin{prop}\label{multi} Let $S'$ be an $R$-scheme and $\chi_1,\ldots, \chi_i: G\to \Gamma(S', \O^*_{S'}) $
characters of $G$. With notations as above, there is a canonical 
homomorphism of invertible sheaves over $Y'$
\begin{equation*}
\mu: \O_{X',  \chi_1}
\otimes_{\O_{Y'}}\cdots\otimes_{\O_{Y'}}\O_{X',\chi_i}\to \O_{X',\chi_1\cdots \chi_i},
\end{equation*}
which is an isomorphism over the complement 
of the branch locus of $\pi': X'\to Y'$. In particular,
the base change $\mu_K$ is an isomorphism over 
the generic fiber $Y_K$.  
\end{prop}

\begin{Proof}
  The homomorphism $\mu$ is given by restricting the $i$-fold multiplication 
\begin{equation*}
\O_{X'}\otimes_{\O_{Y'}}\cdots \otimes_{\O_{Y'}}\O_{X'}\to \O_{X'} \end{equation*}
of the sheaf of rings $\O_{X'}$.  (To simplify the notation we write $\O_{X'}$ instead of
$(\pi')_*(\O_{X'})$.) Over  the complement ${\cal V}'$ of the branch locus of $\pi'$ the morphism
$\pi'_{|{\cal V}'}: {\cal U}'=\pi'^{-1}({\cal V}')\to {\cal V}'$ is a $G$-torsor
and so $\pi'^{-1}({\cal V}')\times_{{\cal V}'}\pi'^{-1}({\cal V}')\simeq G\times \pi'^{-1}({\cal V}')$. A standard argument
using base change by the finite \'etale morphism $\pi'_{|{\cal V}'}$ and descent shows that the restriction of  $\mu$ 
on ${\cal V}'$ is an isomorphism (see for example, [P3, (2.3)]).
 \endproof
\end{Proof}
\medskip

Now consider $I\subset \{1,\ldots, n\}$. Recall that the scheme morphism 
$m_I: (G^D_R)^{n}\to G^D_R$ induced by the algebra homomorphism
$\Delta_{I}$ of \S \ref{cube} is given by $(x_1,\ldots, x_{n})\mapsto \sum_{i\in I}x_i$.
(Recall that we are using additive notation for the group operation on the points of $G^D_R$.
However, such points correspond to characters of $G$ and the group operation 
is actually given by multiplication of characters.)
Also the $i$-th projection morphism $p_{i}: (G^D_R)^{n}\to G^D_R$ is induced by $\Delta_{\{i\}}$.
Now set $\L_I=\bigotimes_{i\in I} p^*_{i}(\L) $. Suppose now that, in the construction of Proposition \ref{multi},
we take $S'=(G^D_R)^n$ and for $i\in I$ we let $\chi_i$ be
the $(G^D_R)^n$-valued character of $G$ which is the ``universal" character $\chi_{\rm u}: G\to R[G]^*$;
$\chi_u(g)=g$ on 
the $i$-th factor and the trivial character on all the other factors.  
Then Proposition \ref{multi} applied to this situation implies:

\begin{cor}\label{multiI}
For each $n\geq 1$ and $I\subset \{1,\ldots, n\}$, there is a canonical homomorphism
\begin{equation*}
\mu^I: \L_I\xrightarrow {\ }  m^*_I(\L)
\end{equation*}
of invertible sheaves over $(G^D_R)^n_Y=Y\times_R(G^D_R)^n$ which is an 
isomorphism over the base change of the complement 
of the branch locus of $\pi: X\to Y$. In particular,
the base change $\mu^I_K$ is an isomorphism over 
the generic fiber $Y_K\times_K(G^D_K)^n$. \endproof
\end{cor}

Proposition \ref{multi} also implies:
\begin{cor}\label{Gpower} 
Let $\chi: G\to \Gamma(S',\O_{S'})$ be a character of $G$. Then there
is a canonical homomorphism of invertible sheaves over $Y'=Y\times_RS'$
\begin{equation*}
\nu: \O_{X',\chi}^{\otimes \#G}\to \O_{X',\chi^{\#G}}=\O_{X',1}=\O_{Y'}
\end{equation*}
which is an isomorphism over the complement 
of the branch locus of $\pi': X'\to Y'$.\endproof
\end{cor}
\smallskip

In the remainder of this paragraph, we assume (in addition to the hypotheses
imposed in the beginning of \S \ref{mainsec}) that $R$ is a complete 
discrete valuation ring with perfect residue field of characteristic prime to $\#G$,
and that $R$ contains a primitive $\#G$-th root of unity $\zeta$.  

Suppose that $x$ is a codimension $1$ 
point of $X$ with $\pi(x)=y$. Our hypothesis (T) implies that the inertia subgroup $I_x$  is trivial
unless $y$ maps to the closed point of $\Spec(R)$.
The action of $I_x$ on the cotangent space ${\mathfrak m}_{X,x}/{\mathfrak m}_{X, x}^2$
defines a faithful character $\phi_x: I_x\to k(x)^*$ with values in the roots of unity of the residue field $k(x)$.
Since $R^*$ contains a root of unity $\zeta$ of order $\#G$ by assumption and $\#G$ is prime
to the characteristic of $k(x)$ we may view $\phi_x$ as taking values in the subgroup $<\zeta>$
generated by $\zeta$.
If $\chi: G\to \ <\zeta>\subset R^*$ is another character of $G$, the restriction of $\chi$ 
to $I_x$ has the form $\phi_x^{n(\chi, x)}$ for a unique integer $n(\chi, x)$ in the range 
$0\leq n(\chi, x)< \#I_x$. Since all points of $X$ over $y$ are conjugate
under the action of $G$, we can see that $I_x\subset G$ and $n(\chi, x)$ depend only $y$.
Often, we will denote them by $I_y$ and $n(\chi, y)$ respectively.
Set
\begin{equation}\label{defig}
g(\chi, y)=-\frac{n(\chi, y)}{\# I_y}.
\end{equation}

\begin{lemma}
Under the above assumptions, the map $\nu: \O_{X,\chi}^{\otimes \#G}\to  \O_{Y}$
identifies $\O_{X,\chi}^{\otimes \#G}$ with $\O_{Y}(F(\chi))$ where
\begin{equation*}
F(\chi)=\sum_{y} \#G\cdot g(\chi, y)\cdot y
\end{equation*}
is the divisor with $y$ running over the finite set of codimension 
$1$ points of $Y$ that are contained in the special fiber of $Y\to \Spec(R)$.
\end{lemma}

\begin{Proof}
Let $y$ be a codimension $1$ 
point of $Y$ which is contained in the special fiber.
It is enough  to prove the statement after replacing $Y$ by an 
\'etale neighborhood of $y\in Y$. Then, by using (\ref{inducedI}), we can
 assume that for every $x\in X$ with $\pi(x)=y$, the decomposition subgroup of $x$
is equal to the inertia subgroup $I_x$.
Suppose that $a$ is a local section of $\O_X$ in a neighborhood of $\pi^{-1}(y)$ such that
$g\cdot a=\chi(g)a$ for $g \in G$, so that $a$ defines a local section of $\O_{X,\chi}\subset \pi_*(\O_X)$ in
a neighborhood of
$y$.   
The subgroup $I_x$ acts on $a$ via the character $\chi_{|I_x}=\phi_x^{n(\chi, x)}$. 
Let   
$\varpi_x$, resp. $\varpi_y$, be uniformizers of the discrete valuation rings $\O_{X,x}$, resp. $\O_{Y,y}$.
Since  $I_x$ acts on the quotient $\varpi_x^j \O_{X,x}/\varpi_x^{j+1} \O_{X,x}$ via
the character $\phi_x^j$ for $j \ge 0$, we see that  $a$ is in $ \varpi_x^{n(\chi,x)} \O_{X,x}$. 
On the other hand, we can choose a local section $a'$ of $\O_{X}$ in a neighborhood of $\pi^{-1}(y)$
which is congruent to $\varpi_x^{n(\chi,x)}$ mod
$\varpi_x^{n(\chi,x)+1} \O_{X,x}$, and which has very high
valuation at every other point of $X$ which lies over $y$.  Define $\alpha = e_\chi \cdot a'$ where
\begin{equation*}
\label{eq:idempotent} e_\chi = \frac{1}{\# G} \sum_{g \in G} \chi(g)^{-1}   g
\end{equation*}
is the idempotent of $\chi$.
Then $\alpha$ satisfies $g\cdot \alpha=\chi(g)\alpha$, and we have
\begin{equation*}
\label{eq:alphacong}
\alpha \equiv \frac{1}{\# G} \sum_{g \in I_x} \chi(g)^{-1}   (g\cdot \varpi_x^{n(\chi,x)}) \equiv 
 \frac{ \# I_x}{\#  G}\, \varpi_x^{n(\chi,x)} 
\, \mathrm{mod}\,
\varpi_x^{n(\chi,x)+1} \O_{X,x} ,
\end{equation*}
since $g\cdot \varpi_x \equiv \phi_x(g) \varpi_x\, {\rm mod}\, \varpi_x^2 \O_{X,x}$ and $\chi(g) = \phi_x(g)^{n(\chi,x)}$
for $g \in I_x$.  By the definition of the map $\nu$, its image 
is the $\O_Y$-ideal sheaf with local sections generated by the  $\#G$-th powers of the local sections
of $\O_{X,\chi}\subset \pi_*(\O_X)$. Since $\varpi_y\O_{X,x}=\varpi_x^{\# I_x}\O_{X,x}$, the above considerations 
now show that the stalk of this ideal
sheaf at $y$ is $\varpi_y^{-g(\chi, y)\cdot \#G}\O_{Y,y}$. 
\endproof
\end{Proof}

\begin{cor}\label{phimuI}
 Let $\phi: G^n\to R^*$ be a character and 
denote by $\phi_i$ the restriction of $\phi$ to the $i$-th factor of $G^n$.
For a subset $I\subset \{1,\ldots, n\}$ denote by $\phi^*(\mu^I)$
the base change of $\mu^I$ of Corollary \ref{multiI}
by the morphism $Y\to Y\times_R(G^D_R)^n=(G^D_R)^n_Y$ induced by base change via $\phi: R[G^n]\to R$. Then   
the $\#G$-th tensor power $\phi^*(\mu^I)^{\otimes \#G}$ of $\phi^*(\mu^I)$  
is identified with the natural injection of invertible sheaves
\begin{equation*}
\O_Y\bigl(\sum_{i\in I} F(\phi_i)\bigr)\to \O_Y\bigl(F(\prod_{i\in I}\phi_i)\bigr).
\end{equation*}
\end{cor}

\subsection{The determinant of cohomology}\label{det}
We continue with the assumptions and notations 
of the beginning of \S \ref{mainsec}. In particular, we assume (T).
Denote by $\ti h: G^D_Y\to G^D_R$ the 
base change of $h$ by $G^D_R\to \Spec(R)$. Recall that we set $\L=\pi_*(\O_X)$ (an invertible sheaf on $G^D_Y$).
If $\Gg$ is a locally free coherent $\O_Y$-module we can 
consider the square of the determinant of cohomology
\begin{equation*}
\delta(\Gg\otimes_{\O_Y}\L)=(\det{\rm R}\ti h_*(\Gg\otimes_{\O_Y}\L))^{\otimes 2};
\end{equation*}
this is a line bundle on $G^D_R$. (If the function
${\rm rank}({\rm R}\ti h_*(\Gg\otimes_{\O_Y}\L))$ is always even on $G^D_R$ we do not have to consider
the square $(\det{\rm R}\ti h_*(\Gg\otimes_{\O_Y}\L))^{\otimes 2}$: In this case, all the 
arguments below can be carried out for $\det{\rm R}\ti h_*(\Gg\otimes_{\O_Y}\L)$
instead of $\delta(\Gg\otimes_{\O_Y}\L)$.)
In order to simplify the exposition we will now 
identify in our notation invertible sheaves (line bundles)  over $G^D_R$ with the corresponding 
locally free rank $1$ $R[G]$-modules of their sections; this should not cause any confusion.
We will also use the same notation $\delta$ to denote the functor
given by the square of the determinant of cohomology over various base schemes without distinction.
Also, in a further attempt to lighten the presentation, we will first concentrate our 
discussion to the case $\Gg=\O_Y$ for simplicity;
the general case, which is similar, will be discussed later.

Let now $s$ be any generator of the $K[G]$-module $\delta(\L)_K$. Our plan is to 
calculate the isomorphism class of the image 
$\Theta_{d+2}((\delta(\L), s))=(\Theta_{d+2}(\delta(\L)), \Theta_{d+2}(s))$ 
in the category ${\cal C}_R(G;d+2)$ in terms of data obtained from the
branch locus of the cover $\pi$. 
Under our assumptions,
the  cover of generic fibers $\pi_K$ is unramified
and hence a $G$-torsor. Hence, it follows from the results of  
[Du] and [P3] (see below) that the $K[G]$-module $\delta(\L)_K$ supports a canonical $n$-cubic
structure $\gamma$ with $n=d+2$. Therefore, there is a distinguished $K[G^{n}]$ generator $\gamma(1)$ of 
$\Theta(\delta(\L)_K)\simeq \Theta(\delta(\L))_K$. (From hereon we will omit the subscript $n=d+2$.) 
We start with the following fundamental observation: 

{\sl It is  enough to calculate 
the isomorphism class of the pair $(\Theta(\delta(\L)), \gamma(1))$ in the group
$C_R(G; n)$ in terms of data obtained from the
branch locus of the cover $\pi$.}

Indeed, since $\gamma$ defines an $n$-cubic structure
on $\Theta(\delta(\L))_K$, if $s$ is a $K[G]$-generator
of $\delta(\L)_K$, then by the definition, the ``ratio" $\Theta(s)/\gamma(1)\in K[G^{n}]^*$ 
is an $n$-cubic element. Hence, the pair $(\Theta(\delta(\L)), \gamma(1))$ is 
isomorphic to $\Theta((\delta(\L), s))=(\Theta(\delta(\L)), \Theta(s))$ 
in the category ${\cal C}_R(G;n)$. 

 Set 
\begin{equation*}
{\cal D}:=\bigotimes_{I\subset\{1,\ldots, n\}} \delta(\L_I)^{(-1)^{n-\#I}}\ .
\end{equation*}
The main result of [Du] (loc. cit. Theorem 4.2) gives a canonical trivialization
\begin{equation*}
b: R[G^{n}]\xrightarrow{\sim} {\cal D}.
\end{equation*}
(This follows from the existence of a $d+2$-cubic structure on the functor
from line bundles on $G^D_Y$ to line bundles on $G^D_R$
given by the square of the determinant of cohomology; see 
also [P3, \S4] especially loc. cit. Definition 4.2 and \S 4.e. In fact, after wrestling with signs, a harder result
about the determinant of the cohomology -- without squaring and without a condition on the rank
of $\R \ti h_*(\L)$ -- is shown
in [Du]. We are not going to use this more complicated result.
When $d=1$, this trivialization can also be obtained directly
using the Deligne-Riemann-Roch theorem and the bilinearity of the Deligne pairing [De].)
Using this we obtain an isomorphism
\begin{equation*}
\Theta(\delta(\L))\otimes {\cal D}^{-1}\xrightarrow {{\rm id}\otimes b^{-1}} \Theta(\delta(\L))\ .
\end{equation*}
Now notice that by using $\mu^I_K$ in Corollary \ref{multiI}
and the functoriality of the determinant of cohomology
we can obtain an isomorphism
\begin{equation}\label{37}
{\cal D}_K \xrightarrow {\sim} \bigotimes_I(\delta((\L_I)_K))^{(-1)^{n-\#I}} \xrightarrow {\sim}
\bigotimes_{I}(\delta(m_I^* (\L)_K))^{(-1)^{n-\#I}}\xrightarrow{\sim}\Theta(\delta(\L)_K)
\end{equation}
where the latter isomorphism comes from the identification
\begin{equation*}
\Theta(\delta(\L))=\bigotimes_Im_I^*(\delta(\L))^{(-1)^{n-\#I}}\simeq 
\bigotimes_I \delta(m^*_I(\L))^{(-1)^{n-\#I}}.
\end{equation*}
Hence, we obtain a trivialization
\begin{equation}\label{38}
K[G^{n}]\xrightarrow {\sim }\Theta(\delta(\L)_K)\otimes {\cal D}^{-1}_K.
\end{equation}
Consider the composition
\begin{equation}\label{compos}
\gamma: K[G^{n}]\xrightarrow {\sim }\Theta(\delta(\L)_K)\otimes {\cal D}^{-1}_K
 \xrightarrow {{\rm id}\otimes b^{-1}_K} \Theta(\delta(\L)_K).
\end{equation}
By [P3, Theorem 4.7] and its proof, the trivialization $\gamma$ 
is an $n$-cubic structure   on $ \delta(\L)_K$.  

Now let us consider the image $E$
of the $R[G^{n}]$-submodule $\Theta(\delta(\L))\otimes 
{\cal D}^{-1}\subset (\Theta(\delta(\L))\otimes {\cal D}^{-1})_K\simeq 
\Theta(\delta(\L)_K)\otimes {\cal D}^{-1}_K$
under the inverse of the isomorphism (\ref{38}) above. 
It follows that $E$ is a locally free $R[G^{n}]$-submodule
of $K[G^{n}]$. There is a commutative diagram
\begin{equation}\label{diag1}
\begin{matrix}
K[G^{n}]&\xrightarrow {\sim }&\Theta(\delta(\L)_K)\otimes {\cal D}^{-1}_K&
 \xrightarrow {{\rm id}\otimes b^{-1}_K}&  \Theta(\delta(\L)_K)\\
\cup&&\uparrow &&  \uparrow\\
E&\xrightarrow  {\sim }&\Theta(\delta(\L))\otimes {\cal D}^{-1}&
 \xrightarrow {  {\rm id}\otimes b^{-1}  }&   \Theta(\delta(\L))\\
\end{matrix}
\end{equation}
where the two vertical arrows are given by the compositions 
$\Theta(\delta(\L))\otimes {\cal D}^{-1}\subset 
(\Theta(\delta(\L))\otimes {\cal D}^{-1})_K\xrightarrow{\sim} \Theta(\delta(\L)_K)\otimes {\cal D}^{-1}_K$
and $\Theta(\delta(\L))\subset \Theta(\delta(\L))_K\xrightarrow{\sim} \Theta(\delta(\L)_K)$.

\begin{lemma}\label{3.6}
Let $s$ be a $K[G]$-generator of the module $\delta(\L)_K$. Then
the element $(E, 1)$ is isomorphic to $\Theta((\delta(\L), s))$
in the category  ${\cal C}_R(G; n)$.
\end{lemma}

\begin{Proof}
An isomorphism is given by the composition of the second row 
of the diagram (\ref{diag1}). As remarked above, $\gamma(1)$ 
is the image of $1$ under the composition of the first row. 
It follows that $(E, 1)$ is isomorphic to $(\Theta(\delta(\L)), \gamma(1))$,
which we have shown earlier to be isomorphic to $\Theta((\delta(\L), s))$
in the category  ${\cal C}_R(G; n)$.\endproof
\end{Proof}

\medskip
>From the definitions, we have canonical isomorphisms
\begin{equation*}
E \xrightarrow{\sim} \Theta(\delta(\L))\otimes{\cal D}^{-1}\xrightarrow{\sim}\bigotimes_{I}\bigl(\delta(m^*_I(\L))\otimes
\delta(\L_I)^{-1}\bigr)^{(-1)^{n-\#I}} =\prod_{I}E_I^{(-1)^{n-\#I}} 
\end{equation*}
when 
\begin{equation*}
E_I=\delta(m^*_I(\L))\otimes \delta(\L_I)^{-1}.
\end{equation*}
Here $E_I$ is a locally free rank $1$ $R[G^n]$-module. We also have a canonical isomorphism
\begin{equation*}
E_I\xrightarrow{\sim}\delta(\L_I\xrightarrow {\mu^I} m_I^*(\L) )
\end{equation*}
in which we think of $\mu^I$ as a perfect complex of $\O_Y[G^n]$-modules
with terms at positions $-1$ and $0$.
Since $\mu^I_K$ is an isomorphism, $\delta(\mu^I_K)$
gives a trivialization $K[G^{n}]\xrightarrow{\sim} (E_I)_K$ that we can use
to identify $E_I$ with a locally free $R[G^{n}]$-submodule of $K[G^{n}]$.
Using these identifications, we obtain
\begin{equation*}
E=\prod_{I}E_I^{(-1)^{n-\#I}}
\end{equation*}
as locally free $R[G^n]$-submodules (``fractional ideals") of $K[G^n]$.
We will obtain our main result by calculating 
$E_I$ and $E$ under some additional assumptions.

\begin{Remarknumb} \label{generalG}
{\rm The above arguments readily extend from the case $\Gg=\O_Y$ to the case of a general 
coherent locally free $\O_Y$-module $\Gg$. Indeed, by [Du, \S 4.7],
the functor from line bundles on $G^D_Y$ to line bundles on $G^D_R$ 
given by $\M\mapsto \delta(\Gg\otimes_{\O_Y}\M)$ has a canonical $n$-cubic structure. 
This provides a canonical trivialization of
\begin{equation}\label{trivG}
{\cal D}(\Gg):=\bigotimes_I\delta(\Gg\otimes_{\O_Y}\L_I)^{(-1)^{n-\#I}}
\end{equation}
over $(G^D_R)^n$.
Now, in the same way as we have seen above, this trivialization combined with the isomorphisms
$\mu^I_K$ provides an $n$-cubic 
structure on the generic fiber $\delta(\Gg\otimes_{\O_Y}\L)_K$. The rest of the argument is 
also identical:
In the end we obtain locally free $R[G^n]$-submodules $E_I(\Gg)$, $E(\Gg)\subset K[G^n]$
such that $(E(\Gg), 1)$ is isomorphic to the image $\Theta((\delta(\Gg\otimes_{\O_Y}\L), s))$
in ${\cal C}_R(G;d+2)$.}
\end{Remarknumb}

\subsection{Branch divisors and Riemann-Roch} \label{branch}
We continue with the assumptions and notations given at the beginning of \S \ref{mainsec}. 
However, in this paragraph, we will assume in addition that $R$ is a complete discrete valuation ring 
which has perfect residue field $k$ of characteristic prime to $\#G$ and which contains a primitive root of unity $\zeta$
of order equal to $\#G$. Recall also $n=d+2$. 

Let $\phi: G^n\to\ <\zeta>\subset R^*$ be a character and denote by $\phi_i$ the restriction
of $\phi$ to the $i$-th factor of $G^n$. Recall that for $I\subset \{1,\ldots, n\}$ we denote by 
$\phi^*(\mu^I)$ the base change of $\mu^I$ by $Y\to (G^D_Y)^n$ given by using $\phi: R[G^n]\to R$. The morphism 
$\phi^*(\mu^I)$ then identifies 
with the ``multiplication"
\begin{equation*}
\bigotimes_{i\in I}\O_{X,\phi_i}\xrightarrow{ }  \O_{X,\prod_{i\in I}\phi_i}
\end{equation*}
which was considered in Proposition \ref{multi}. By functoriality, we have
\begin{equation*}
\phi(E_I)=\delta(\phi^*(\mu^I))\subset K, \quad \phi(E)=\prod_{I}\delta(\phi^*(\mu^I))^{(-1)^{n-\#I}}\subset K .
\end{equation*}
We will calculate the $R$-fractional ideals $\phi(E_I)$ and $\phi(E)$ in $K$
using a localized Riemann-Roch theorem for the morphism $h: Y\to \Spec(R)$.
Denote by $h_s: Y_s\to \Spec(k)$ the special fiber of $h$. Let $\F_{I,\phi}$ be the cokernel 
of $\phi^*(\mu^I)$; it is a coherent sheaf 
of $\O_{Y}$-modules which is supported on the special fiber $Y_s$.
Therefore, it gives a class $[\F_{I,\phi}]$ in the Grothendieck group ${\rm G}_0(Y_s)$ of 
coherent $\O_{Y_s}$-modules on $Y_s$. The following can be deduced from the definition
and basic properties of the determinant of cohomology. (Note that the functor $\delta$ 
in this statement is given by the {\sl square} 
of the determinant of cohomology; this explains the appearance of the factor $2$ in the exponent below.)

\begin{lemma}\label{valua}
Denote by $\varpi$ a uniformizer of $R$. We have
\begin{equation*}
\phi(E_I)=\varpi^{-2\cdot\chi(\F_{I,\phi})}\cdot R,
\end{equation*}
where $\chi(\F_{I,\phi})=(h_s)_*([\F_{I,\phi}])\in {\rm G}_0(\Spec(k))=\Z$ is 
the Euler characteristic of $[\F_{I,\phi}]$.\endproof
\end{lemma}

In what follows we borrow heavily from [Fu]. We refer the reader to loc. cit. for notations and more details.
Suppose that $D_1,\ldots, D_q$ are effective divisors on $Y$ 
with supports $|D_1|,\ldots , |D_q|$
contained in $Y_s$. Then
we can define, following [Fu, \S 2 and \S 17], a bivariant class $[D_1]\cap\cdots \cap [D_q]\in
A^q(Y_s\to Y)$ which induces homomorphisms $A_{*}(Y)\to A_{*-q}(Y_s)$ denoted by $x\mapsto ([D_1]\cap\cdots
\cap [D_q])\cap x$. These homomorphisms are given by the composition  
$$
A_*(Y)\to A_{*-q}(|D_1|\cap \cdots \cap |D_q|)\to A_{*-q}(Y_s).
$$
Here, the first map is defined by [Fu, \S 2] and the second is obtained from the inclusion of the set-theoretic
intersection of supports $|D_1|\cap \cdots \cap |D_q|$ in $Y_s$. The bivariant class $[D_1]\cap\cdots \cap [D_q]$ is 
multiadditive in the $D_i$ and independent of their order ([Fu, \S 2.3, \S 2.4]). (If $z=[Z]\in A_k(Y)$ with $Z$
integral of Krull dimension $k$ then $[D]\cap z$ is the class of the $(k-1)$-cycle on $|D|\subset Y_s$ obtained by restricting the line bundle $\O_Y(D)$ on $Z$ and taking a corresponding divisor 
supported on $|D|\cap Z$. The general definition follows from this by using an inductive procedure.) 
If $D_1=\cdots =D_q=D$, we will denote this class by $[D]^q$.

Now let $\E^\bullet$ be a finite complex of locally free $\O_Y$-modules which is exact off $Y_s$.
Denote by ${\rm ch}_{Y_s}^Y(\E^\bullet)$ the localized Chern character 
of  $\E^\bullet$  defined by the MacPherson graph 
construction following [Fu, \S 18]. (Strictly speaking, the reference [Fu] only covers 
schemes over a base field; the extension to the situation we want is described 
in [GS] or [Ro2, \S 11], see also [Fu, \S 20].)  
By definition, this is a bivariant
class in $A^*(Y_s\to Y)_\Q$ which in particular induces
\begin{equation*}
{\rm ch}_{Y_s}^Y(\E^\bullet)\cap \cdot : A_*(Y)_\Q\to A_*(Y_s)_\Q\ .
\end{equation*}
Suppose that $c$ is a complex of the form 
$\L\to \M$ (at degrees $-1$ and $0$) with $\L$, $\M$  two line bundles on $Y$
which is exact off $Y_s$ and set
\begin{equation*}
{\rm ch}_{Y_s}^Y(c)=
\sum_{q=0}^{d+1} 
{\rm ch}_{Y_s}^{Y, q}(c)\in A^*(Y_s\to Y)_\Q=\oplus_{q=0}^{d+1}A^q(Y_s\to Y)_\Q
\end{equation*}
We have (see [Ro2, Theorem 11.4.6, Theorem 12.3.1]):

\begin{prop}\label{chern}
a) ${\rm ch}_{Y_s}^{Y, 0}(c)=0$.

b) Suppose that $c$ is the inclusion $\O_Y\to \O_Y(D)$ for $D$ an effective divisor
with support $|D|\subset Y_s$. Then 
\begin{equation*}
{\rm ch}^{Y}_{Y_s}(c)=\sum_{q=1}^{d+1} \frac{[D]^q}{q!}.
\end{equation*}

c) Similarly, if $c$ is the inclusion $\O_Y(-D)\to \O_Y$ for $D$ an effective divisor
with support $|D|\subset Y_s$, then 
\begin{equation*}
{\rm ch}^{Y}_{Y_s}(c)=-\sum_{q=1}^{d+1} \frac{(-1)^q[D]^q}{q!}.
\end{equation*}

d) If $c^m: {\L^{\otimes m}}\to {\M^{\otimes m}}$ is the $m$-fold tensor product of $c$,
then 
\begin{equation*}
{\rm ch}_{Y_s}^{Y, q}(c^m)=m^q {\rm ch}_{Y_s}^{Y, q}(c).
\end{equation*}

e) If $\Gg$ is a locally free coherent $\O_Y$-module, then ${\rm ch}^{Y}_{Y_s}(\Gg\otimes_{\O_Y}c)
={\rm ch}(\Gg)\cap {\rm ch}^{Y}_{Y_s}(c)$, where ${\rm ch}(\Gg)\in A^*(Y)_\Q$ is the usual Chern character
of $\Gg$.\endproof

\end{prop}

Now let us return to the calculation of the fractional ideals $\phi(E_I)$, $\phi(E)$.
Since $Y$ is assumed regular and $h$ projective,   an embedding $i: Y\to {\bf P}^m_R$
produces an exact sequence of coherent $\O_Y$-modules
\begin{equation*}
0\to N^\vee_{Y|{\bf P}^m_R}\to i^*\Omega^1_{{\bf P}^m_R{/R}}\to \Omega^1_{Y/R}\to 0
\end{equation*}
with $N^\vee_{Y|{\bf P}^m_R}$, $i^*\Omega^1_{{\bf P}^m_R/R}$ locally free.
(Here $N^\vee_{Y|{\bf P}^m_R}$ is the dual of the normal bundle $N_{Y|{\bf P}^m_R}$ 
of the embedding.)
Denote by ${\rm Td}(h)={\rm td}(h)\cap [Y]$ the Todd class in $A_*(Y)_\Q$
of the relative  tangent complex $[(i^*\Omega^1_{{\bf P}^m_R/R})^\vee\to N_{Y|{\bf P}^m_R}]$
(at degrees $0$ and $1$).

\smallskip
\begin{thm}\label{chtd}
With the above assumptions and notations, consider the function 
$T_\pi$ on $R$-valued $1$-dimensional characters of $G$
given by
\begin{equation*}
\psi\mapsto  T_\pi(\psi):=\sum_{q=1}^{d+1}\sum_{(y_1,\ldots, y_q)}\frac{\prod_{j=1}^qg(\psi, y_j)}
{q!}{\rm deg}\Big(\big(\cap_{j=1}^q[ y_j]\big)\cap {\rm Td}_{q}(h)\Big)\in \Q,
\end{equation*}
where we write ${\rm deg}$ for the degree of zero cycles over $k$
and $(y_1,\ldots, y_q)$ runs over all $q$-tuples of codimension $1$ points of $Y$ that lie in $Y_s$.
(We use the same symbol for the corresponding effective divisors.)
Extend linearly the function $T_\pi$  to the ring of $R$-valued characters of $G$, i.e to integral linear 
combinations of $R$-valued $1$-dimensional characters.
Then we have
\begin{equation*}
\phi(E)=\varpi^{-2\cdot T_\pi(\Theta^D(\phi))}\cdot R\subset K
\end{equation*}
where for $\phi=\phi_1\otimes\cdots\otimes \phi_{d+2}$ a $1$-dimensional character of $G^{d+2}$ we have
\begin{equation*}
\Theta^D(\phi)=(\phi_1-1)\cdots (\phi_{d+2}-1)=\sum_{I}(-1)^{d+2-\#I}\prod_{i\in I} \phi_i\ .
\end{equation*}
\end{thm}

\begin{Proof} We will see that this follows from Lemma \ref{valua} and:

\begin{thm} \label{RR}(Riemann-Roch theorem)
\begin{equation*}
\chi(\F_{I,\phi})=(h_s)_*\left(({\rm ch}^Y_{Y_s}(\phi^*(\mu^I))\cap {\rm Td}(h))_{0}\right).
\end{equation*}
\end{thm}

In this equality, the map $(h_s)_*$ on the right hand side is the 
push forward of zero-cycles $A_0(Y_s)_\Q$ 
to $A_0(\Spec(k))_\Q=\Q$; i.e given by the degree of zero cycles over $k$.
This special case of a ``localized" Riemann-Roch theorem follows from [Ro2, Theorem 12.5.1] and [Ro2, Theorem 12.6.1].
It can also be derived following the proof of [Fu, Theorem 18.2 (1)] by considering the morphism 
$h_s: Y_s\to \Spec(k)$ as a morphism of schemes over $S=\Spec(R)$. 
(This latter reference gives a similar result for schemes over a base 
$S$ which is a non-singular scheme over a 
field.)
\smallskip

Now, let us deduce Theorem \ref{chtd} from Theorem \ref{RR} by calculating 
the $0$-th component of $\sum_{I}(-1)^{n-\# I}  {\rm ch}^Y_{Y_s}(\phi^*(\mu^I))\cap {\rm Td}(h)  $: 
Using Proposition \ref{chern} (a) and (d), we can write
\begin{equation}\label{Gtod+1}
({\rm ch}^Y_{Y_s }(\phi^*(\mu^I))\cap {\rm Td} (h))_{0}=
\sum_{q=1}^{d+1} {\rm ch}_{Y_s}^{Y, q}(\phi^*(\mu^I))\cap {\rm Td}_{q}(h)=
\end{equation}
\begin{equation*}
=\sum_{q=1}^{d+1} \frac{{\rm ch}_{Y_s}^{Y, q}(\phi^*(\mu^I)^{\otimes\#G})}{(\#G)^q}\cap {\rm Td}_{q}(h).
\end{equation*}
Corollary \ref{phimuI} now identifies the complex $\phi^*(\mu^I)^{\otimes\# G}$ with 
\begin{equation*}
c'_I: \O_Y(\sum_{i\in I} F(\phi_i))\to \O_Y(F(\prod_{i\in I}\phi_i)).
\end{equation*}

Let us  write $D_I=-F(\prod_{i\in I}\phi_i)$, $D_I+D'_I=-\sum_{i\in I}F(\phi_i)$
(these are both effective divisors supported on the special fiber $Y_s$). There is an exact sequence of complexes:
$$
\begin{matrix}
0\to&\O_Y(-D_I) & \xrightarrow{}  &\O_Y(-D_I)\oplus \O_Y &\xrightarrow {} & \O_Y &\to  0\\
&&&&&&\\
&c'_I\ \uparrow\ \  &&\uparrow {\rm id}\oplus i''& & \uparrow i\\
&&&&&&\\
0\to&\O_Y(-D_I-D_I') & \xrightarrow{}&  \O_Y(-D_I)\oplus \O_Y(-D_I-D_I') &\xrightarrow {} &\O_Y(-D_I) &\to  0\\
\end{matrix}
$$
with $c'_I$, $i$, $i''$ the natural injective homomorphisms.
Therefore, by [Fu, Proposition 18.1], we have:
\begin{equation}\label{addi}
{\rm ch}^Y_{Y_s}(c'_I)+{\rm ch}^Y_{Y_s}(\O_Y(-D_I)\xrightarrow{i} 
\O_Y)={\rm ch}^Y_{Y_s}(\O_Y(-D_I-D'_I)\xrightarrow{i''} \O_Y).
\end{equation}
Recall $n=d+2$.
Applying Proposition \ref{chern} (c) and telescoping using the identity
(for given $q<n$)
\begin{equation*}
\sum_{I\subset\{1,\ldots, n\}}(-1)^{\#I}\left(\sum_{i\in I}X_i\right)^q=0
\end{equation*}
gives
\begin{equation*}
\sum_{I\subset \{1,\ldots, n\}}(-1)^{n-\# I}{\rm ch}^Y_{Y_s}(\O_Y(\sum_{i\in I}F(\phi_i))\to \O_Y )=0.
\end{equation*}
This translates to 
\begin{equation*}
\sum_{I\subset \{1,\ldots, n\}}(-1)^{n-\# I}{\rm ch}^Y_{Y_s}(\O_Y(-D_I-D'_I)\xrightarrow{i''} \O_Y)=0.
\end{equation*}
Hence, using (\ref{addi}) and Proposition \ref{chern} (c), we obtain
\begin{equation*}
\sum_{I}(-1)^{n-\# I}{\rm ch}^Y_{Y_s}(c'_I)=-\sum_{I}(-1)^{n-\# I}{\rm ch}^Y_{Y_s}(\O_Y(-D_I)\to \O_Y))=
\end{equation*}
\begin{equation*}
=\sum_{I}(-1)^{n-\# I}\sum_{q=1}^{d+1}\frac{(-1)^q[D_I]^q}{q!}.
\end{equation*}
 This gives
\begin{equation*}
\sum_{I}(-1)^{n-\# I}{\rm ch}_{Y_s}^{Y, q}(c'_I)=\sum_{I}(-1)^{n-\# I+q}\frac{ [D_I]^q}{q!}.
\end{equation*}
Therefore, using Proposition \ref{chern} (d) and Corollary \ref{phimuI}, we now obtain
\begin{equation*}
\sum_{I}(-1)^{n-\# I}{\rm ch}_{Y_s}^{Y, q}(\phi^*(\mu^I))=
\sum_{I}(-1)^{n-\# I}\frac{{\rm ch}_{Y_s}^{Y, q}(\phi^*(\mu^I)^{\otimes\#G})}{(\#G)^q}=
\end{equation*}
\begin{equation*}
=\sum_{I}(-1)^{n-\# I+q}
\frac{ [D_I]^q}{(\#G)^q q!}.
\end{equation*}
Since by definition $D_I=-F(\prod_i\phi_i)$, by using Lemma 3.5 we find that this is equal to
\begin{equation*}
\sum_{I}(-1)^{n-\# I+q} \sum_{(y_1,\ldots, y_q)}
\frac{{(-\#G)}^q  g(\prod_{i\in I}\phi_i, y_1)\cdots g(\prod_{i\in I}\phi_i, y_q)}
{{(\#G)}^q q!}[y_1]\cap \cdots \cap [y_q]=
\end{equation*}
\begin{equation*}
=\sum_{I}(-1)^{n-\# I}\sum_{(y_1,\ldots, y_q)}
\frac{g(\prod_{i\in I}\phi_i, y_1)\cdots g(\prod_{i\in I}\phi_i, y_q)}
{q!}[y_1]\cap \cdots \cap [y_q]
\end{equation*}
where $(y_1,\ldots, y_q)$ runs over all $q$-tuples of codimension $1$ fibral points of $Y$.
Altogether we obtain that the $0$-th component of $\sum_{I}(-1)^{n-\# I}{\rm ch}^Y_{Y_s}(\phi^*(\mu^I))\cap {\rm Td}(h)$
is
\begin{equation*}\label{formula}
\sum_{I}(-1)^{n-\# I}\sum_{q=1}^{d+1}\sum_{(y_1,\ldots, y_q)}\frac{\prod_{j=1}^qg(\prod_{i\in I}\phi_i, y_j)}
{q!}\left(\cap_{j=1}^q[y_j]\right)\cap {\rm Td}_{q}(h).
\end{equation*}
Theorem \ref{chtd} now follows from Theorem \ref{RR} and Lemma \ref{valua}.
\endproof
\end{Proof}

\subsection{}\label{remG3.10} The results of \S \ref{branch} can  be  extended easily
from the case $\Gg=\O_Y$ to the case of a general 
coherent locally free $\O_Y$-sheaf $\Gg$: We want to calculate
\begin{equation*}
\phi(E(\Gg))=\bigotimes_I\delta\bigl(\phi^*\bigl(\Gg\otimes_{\O_Y}\L_I
\xrightarrow{\Gg\otimes_{\O_Y}\mu^I}\Gg\otimes_{\O_Y}m^*_I(\L)\bigr)\bigr)^{(-1)^{n-\#I}}.
\end{equation*}
The exact same proof now works (one has to also use Proposition \ref{chern} (e))
to obtain that 
\begin{equation}\label{gen3.9}
\phi(E(\Gg))=\varpi^{-2\cdot T_{\pi,\Gg}(\Theta^D(\phi))}\cdot R\subset K,
\end{equation}
where $T_{\pi, \Gg}$ is the $\Q$-valued function on the $R$-valued $1$-dimensional characters 
of $G$ given by 
\begin{equation} \label{genT}
T_{\pi,\Gg}(\psi)=\sum_{l=1}^{d+1}
\sum_{(y_1,\ldots, y_{l})}\frac{\prod_{j=1}^{l}g(\psi, y_j)}{l!} 
\sum^{d+1-l}_{t=0}{\rm deg}\Big({\rm ch}^t(\Gg)\cap 
\big(\cap_{j=1}^{l}[ y_j]\big)\cap {\rm Td}_{t+l} (h)\Big).
\end{equation}
The notation and assumptions here are as in Theorem \ref{chtd}.
\bigskip

\subsection{} 
Here we assume in addition that $X$ and $Y$ are relative 
curves over $\Spec(R)$ ($d=1$). Notice that we have
\begin{multline}
\ \ \ \ \ \ \  {\rm Td}_1(h)=\frac{c_1((i^*\Omega^1_{{\bf P}^m_R/R})^\vee)}{2}-\frac{c_1(N_{Y|{\bf P}^m_R})}{2}=\\
=\frac{c_1(\det(i^*\Omega^1_{{\bf P}^m_R/R})^{-1}\otimes \det(N^\vee_{Y|{\bf P}^m_R}) )}{2}=
-\frac{c_1(\omega)}{2}\ \ \ \ \  \ \ 
\end{multline}
where $\omega\simeq \det(i^*\Omega^1_{{\bf P}^m_R/R})\otimes \det(N^\vee_{Y|{\bf P}^m_R})^{-1}$
is the canonical sheaf for $Y\to \Spec(R)$. Hence, in this case, the formula (\ref{genT})
specializes to 
\begin{multline}\label{surf1}
\ \ \ \ \ \ \ \   T_{\pi, \Gg}(\psi)={\rm rank(\Gg)}\sum_{(y_1,y_2)}\frac{g(\psi, y_1)\cdot g(\psi,y_2)}{2}\deg([y_1]\cap [y_2])-\\
\ \ \ \   -\sum_{y}g(\psi, y)\deg\left({\rm rank(\Gg)}[y]\cap  \frac{c_1(\omega)}{2}-c_1(\Gg)\cap [y] \right).\ \ \ \ \ \ 
\end{multline}
Actually, in this case we have $\deg([y_1]\cap [y_2])=\deg(\O_Y(y_1)_{|y_2})$ 
(the degree of the line bundle $\O_Y(y_1)$ restricted on $y_2$ over $\Spec(k)$). For simplicity, we will
denote this intersection number by $y_1\cdot y_2$ (see [La2, III]). Also, 
using adjunction ([La2, IV \S 4]), we find that $\deg([y]\cap c_1(\omega)) =-\deg([y]\cap[y])+\deg(c_1(\omega_{y/k}))=-y\cdot y-2\chi(y,\O_y)$
(here $\omega_{y/k}$ is the canonical sheaf of $y$ over $\Spec(k)$, see loc. cit.). Then formula (\ref{surf1}) for $\Gg=\O_Y$ can be written:
\begin{equation}\label{surf2}
 T_{\pi }(\psi)= \sum_{(y_1,y_2)}\frac{g(\psi, y_1)\cdot g(\psi,y_2)}{2}(y_1\cdot y_2)+
\sum_{y}\frac{g(\psi, y)}{2}\left(y\cdot y+2\chi(y,\O_y)\right).  
\end{equation}

\subsection{The main theorem}\label{mainthmpar} Let $R$ be the ring of integers of a number field $K$. 
We continue with the assumptions and notations 
of the beginning of \S \ref{mainsec}. In fact, in addition to (T) we now also assume:
\smallskip

(D) the residue field characteristic of each  point 
of $Y$ which ramifies in $\pi: X\to Y$ is 
relatively prime to the order of the group $G$.
\smallskip

Denote by $S$ the finite set of rational primes such that the cover $\pi: X\to Y$
is only ramified at points above $S$. By our assumption (D), $p\in S$ implies $p{\not| \#G}$.
Denote by $S_K$ the set of places of $K$ that lie above $S$.

Suppose $\Gg$ is a locally free coherent $\O_Y$-sheaf on $Y$ and
 consider the $G$-sheaf $\F=\pi^*\Gg$ on $X$. Recall the homomorphism
\begin{equation*}
\Theta=\Theta_{d+2}: \Cl(R[G])={\rm Pic}(R[G])\to C_R(G;d+2)
\end{equation*}
of \S \ref{cubecat}. For a finite place $v$ of $K$, we denote by $\varpi_v$ a uniformizer of 
the completion $R_v$ and fix an algebraic closure $\bar K_v$ of its fraction field $K_v$.
Recall that, as in \S \ref{ideles}, any finite id\`ele $(a_v)_v\in {\bf A}^*_{f, K[G^{d+2}]}$
gives the element $(\cap_v( R_v[G^{d+2}]a_v\cap K[G^{d+2}]), 1)$ of $C_R(G;d+2)$.  
Now let $v\in S_K$ (then $(v,\#G)=1$) and denote by $R'_v$ the 
complete discrete
 valuation ring $R'_v\subset \bar K_v$ obtained by adjoining to $R_v$ a primitive root of unity of order 
equal to $\#G$. Then $\varpi_v$ is also a uniformizer for $R'_v$.
Let us consider the cover $\pi'_v: X\otimes_R {R'_v}\to Y\otimes_R {R'_v}$ 
obtained from $\pi$ by base change. Since $R'_v$ has residue field 
characteristic prime to $\#G$ and contains a primitive $\#G$-th root of unity
we can now apply the constructions and results of 
paragraphs  \ref{branch} and \ref{remG3.10} to the cover $\pi'_v$ and the sheaf $\Gg\otimes_RR'_v$. For simplicity, we will denote by
$T_{v,\Gg}$ the function $T_{\pi'_v, \Gg\otimes_RR'_v}: {\rm Ch}(G)_v\to \Q$ 
given by (\ref{genT}). Recall the isomorphism 
\begin{equation}\label{hom}
K_v[G^{d+2}]^*\xrightarrow {\sim} {\rm Hom}_{{\rm Gal}(\bar K_v/K_v)}({\rm Ch}(G^{d+2})_v, \bar K_v^*)
\end{equation}
given by evaluating characters of $G^{d+2}$.  Also recall that if 
$\phi=\phi_1\otimes\cdots \otimes \phi_{d+2}$ is a character of $G^{d+2}$
given by a $d+2$-tuple $(\phi_i)_i$ of $1$-dimensional $\bar K_v$-valued characters of $G$,
we have $\Theta^D(\phi)=(\phi_1-1)\cdots (\phi_{d+2}-1)\in {\rm Ch}(G^{d+2})_v$.
Let us now observe that the definition of the function $T_{v,\Gg}$ implies that, 
for all $\sigma\in {{\rm Gal}(\bar K_v/K_v)}$ 
and $\psi\in {\rm Ch}(G)_v$, we have
$$
T_{v,\Gg}( \psi)=T_{v,\Gg}(\psi^\sigma).
$$
Indeed, this can be derived from the equality $g(\psi, y_j)=g(\psi^\sigma, y_j^\sigma)$
in which $y_j^\sigma$ is the image of the irreducible component $y_j$ under the action of 
$\sigma\in {\rm Gal}(\bar K_v/K_v)$. This, in turn, follows from the identity
\begin{equation*}
 \left(\frac{g(\varpi)}{\varpi}\right)^\sigma=
\frac{g( \varpi^\sigma)}{ \varpi^\sigma}\ , 
\end{equation*}
with $g\in I_{y_j}$ and  $\varpi$, $\varpi^\sigma$ uniformizers of the local rings at $y_j$ and $y_j^\sigma$ respectively,
which is true since the $G$-action is ``defined over $R$".
Hence, the map $\phi\mapsto \varpi_v^{-T_{v,\Gg} (\Theta^D(\phi))}$ gives a function in 
${\rm Hom}_{{\rm Gal}(\bar K_v/K_v)}({\rm Ch}(G^{d+2})_v, \bar K_v^*)$.

\begin{thm}\label{mainthm} With the above assumptions and notations,
\begin{equation*}
\Theta(2\cdot\bar\chi^P(X,\F))=(\cap_v( R_v[G^{d+2}]\lambda_v\cap K[G^{d+2}]), 1),
\end{equation*}
where $(\lambda_v)_v\in 
{\bf A}^*_{f, K[G^{d+2}]}$
is  the (unique) finite id\`ele which is such that
\begin{equation}\label{345}
\phi(\lambda_v)=\begin{cases}
1, &\text{if $v \not\in S_K$ ;}\\
\varpi_v^{-2\cdot T_{v,\Gg} (\Theta^D(\phi))}, &\text{if $v\in S_K$\ \ ,}
\end{cases}
\end{equation}
for all $\bar K_v$-valued characters $\phi $ of $G^{d+2}$. 

If the
``usual" Euler characteristic  $\chi(Y,\Gg)=\sum_i(-1)^i{\rm rank}_R[{\rm H}^i(Y,\Gg)] $
is {\sl even}, then we can eliminate both occurrences of the factor $2$ from the
statement: $\Theta(\bar\chi^P(X,\F))$ is then given by the id\`ele $(\lambda'_v)_v$
with $\phi(\lambda'_v)=1$ if $v\not\in S_K$, and $\phi(\lambda'_v)=\varpi_v^{ -T_{v,\Gg} (\Theta^D(\phi))}$
if $v\in S_K$.
\end{thm}

\begin{Proof}
By (\ref{identity}) it is enough to consider the image of the class of the
determinant of cohomology 
$\det \R\ti h_*(\pi_*(\F))=\det \R\ti h_*(\Gg\otimes_{\O_Y}\L)$ (or of its square) under $\Theta$.
Notice that under the assumption that $\chi(Y,\Gg)=\sum_i(-1)^i{\rm rank}_R[{\rm H}^i(Y,\Gg)]$
is even, the function given by ${\rm rank}(\R\ti h_*(\Gg\otimes_{\O_Y}\L))$ is always even on $G^D_R$.
Indeed, since this function is Zariski locally constant and $G^D_R$ is connected, it
is enough to check this over the generic fiber $\Spec(K)$; there it follows from the Grothedieck-Riemann-Roch theorem
and the fact that $\L_K$ is a torsion line bundle (cf. Cor. \ref{Gpower}) and hence it has
trivial Chern character. Hence in this  ``even" case we do not have to consider
the square.

By Remark \ref{generalG} (cf. Lemma \ref{3.6}), it is now enough to show that, for each finite place $v$ of $K$,
$
E(\Gg)R_v=R_v[G^{d+2}]\lambda_v
$
in  $K_v[G^{d+2}]$. When $v{\not\in}S_K$, $E(\Gg)R_v=R_v[G^{d+2}]$  
by Corollary \ref{multiI} and the definition of $E(\Gg)$. Hence, since $\lambda_v=1$,
$E(\Gg)R_v=\lambda_vR_v[G^{d+2}]$.
When $v\in S_K$, $(v,\# G)=1$, and
we have an isomorphism
\begin{equation*}
R_v[G^{d+2}]\xrightarrow{\sim} \oplus_{\phi}R_v(\phi),\quad a\mapsto (\phi(a))_{\phi},
\end{equation*}
with $\phi$ ranging over the ${\rm Gal}(\bar K_v/K_v)$-orbits of $\bar K_v$-valued characters
of $G^{d+2}$ and $R_v(\phi)\subset R'_v\subset \bar K_v$ the 
complete discrete valuation ring generated over $R_v$ by the values 
of $\phi$. Therefore, to show the statement for $v\in S_K$, it is enough 
to verify that $\phi(E(\Gg))R_v(\phi)=\phi(\lambda_v)R_v(\phi)$.
For that is enough to check
$\phi(E(\Gg))R'_v=\phi(\lambda_v)R_v'$. This now follows from   
Remark \ref{remG3.10} (cf. Theorem \ref{chtd}), and the fact that
the formation of the determinant of cohomology (and therefore also of the 
ideals $E(\Gg)$) commutes with base change.\endproof 
\end{Proof}

  \begin{Remarknumb}
\label{kummer}
{\rm Theorem \ref{mainthm} determines $2\cdot\bar\chi^P(X,\F)$ up to an element of 
${\rm Ker}(\Theta_{d+2})$ from numerical information.  The group ${\rm Ker}(\Theta_{d+2})$ is bounded above in Theorem
\ref{2.7}; in particular, it is trivial if $R = \Z$ and $ d \le 3$.    Note, 
though, that differs from the usual approach in Galois module structure
theory of specifying an id\`ele $\beta$ in ${\bf A}^*_{f, K[G]}$  which represents $2\cdot\bar\chi^P(X,\F)$ 
(e.g. modulo ${\rm Ker}(\Theta_{d+2})$ or a larger subgroup of $\Cl(R[G])$).  In the modular form example 
of the next section
we produce such a $\beta$   
under some additional hypotheses. In general, we can produce an id\`ele in 
${\bf A}^*_{f, K[G]}$ giving the class $2(\#G)^{d+1}\chi^P(X,\F)$ modulo 
${\rm Ker}(\Theta_{d+2})$; see  \S \ref{outputAB} below.}
\end{Remarknumb}

\begin{Remarknumb}\label{Gsheaf}
{\rm Suppose that $\F$ is a general locally free $G$-sheaf on $X$
(i.e a locally free coherent $\O_X$-module with $G$ action compatible
with the action of $X$) which is not necessarily of the
form $\pi^*\Gg$ for an $\O_Y$-sheaf $\Gg$. In this remark we explain how our methods 
can also lead to a calculation of $\bar\chi^P(X,\F)$. Set $\Hh=(\pi_*(\F))^G$. Under our assumptions, this
is a locally free $\O_Y$-sheaf and there is an exact sequence of 
$G$-sheaves on $X$
\begin{equation}\label{piHQ}
0\to \pi^*(\Hh)\to \F\to {{\mathcal Q}}\to 0
\end{equation}
with $\mathcal Q$ supported on the (fibral) branch locus of $\pi$. We have
$\bar\chi^P(X,\F)=\bar\chi^P(X,\pi^*(\Hh))+\bar\chi^P(X,{\mathcal Q})$. The class
$\bar\chi^P(X,\pi^*(\Hh))$ can be calculated using Theorem \ref{mainthm}.
On the other hand, the class $\bar\chi^P(X,{{\mathcal Q}})$ naturally lifts to a class
in the torsion free group $\oplus_{v}{\rm K}_0(k(v)[G])$ which can be calculated (at least when $X$ is 
regular) using the Lefschetz-Riemann-Roch theorem of [BFQ]
(see [CEPT2] for an example of such a calculation). 
}
\end{Remarknumb}

\subsection{}\label{outputAB}
 Let us discuss here how we can deduce Theorem \ref{thm:inputoutput}
of the introduction. Observe that under the assumptions of Theorem \ref{thm:inputoutput} 
$\pi$ is flat by Remark \ref{abhya} (a) and we can write the short exact
sequence (\ref{piHQ}) with ${\mathcal Q}$ supported on the fixed point locus $X'$.
Moreover, we can see that ${\mathcal Q}$ depends only on the pair $(\hat X', \F{|\hat X'})$.
As above, we have $\bar\chi^P(X,\F)=\bar\chi^P(X,\pi^*(\Hh))+\bar\chi^P(X,{\mathcal Q})$.
The class $\bar\chi^P(X,{{\mathcal Q}})$ is represented by a finite module
of finite projective dimension of order supported on the image of $X'$ 
in $\Spec(\Z)$. Hence, we can reduce the proof to the case
 that $\F$ is of the form $\pi^*\Gg$
for $\Gg$ a locally free coherent sheaf on $Y$. Then (a) follows by combining
Theorem \ref{mainthm} with Theorem \ref{2.7}.
To show (b) observe that the character
function $\psi\mapsto (\#G)^{d+1}T_{\pi,\Gg}(\psi)$
takes values in $\Z$. (This follows from the 
definition of $T_{\pi,\Gg}(\psi)$ and the localized Riemann-Roch
theorem.)
Therefore, we can use this character function to define by a formula as in 
(\ref{345}) an id\`ele in ${\bf A}^*_{f, \Q[G]}$ 
whose image under the map $\Theta_{d+2}$ is $2(\#G)^{d+1}\cdot\bar\chi^P(X,\F)$.
By construction, this id\`ele represents a finite module of finite projective 
dimension whose order is supported on the primes below the branch locus
and whose class in $\Cl(\Z[G])$ agrees with $2(\#G)^{d+1}\cdot\bar\chi^P(X,\F)$
modulo ${\rm ker}(\Theta_{d+2})$. The result now follows from Theorem \ref{2.7}.

\bigskip
\bigskip

\section{Galois structure of  modular forms}\label{forms}
\label{s:modforms}

The object of this section is to prove Theorem \ref{qua1}.  With the assumptions of
that Theorem, 
let ${\mathfrak M}_{\Gamma_1(p)}$ be the moduli stack classifying triples 
$(E, C, \gamma)$ where $E\to S$ is a generalized elliptic curve (see [DR]), $C$
a locally free rank $p$ subgroup of the smooth locus $E^{\rm sm}$ and $\gamma: (\Z/p\Z)_S\to C^D$
a ``generator" of the Cartier dual of $C$;  we ask that $C$ intersects every  irreducible component
of every geometric fiber of $E\to S$. (Here generator is meant
in the sense of [KM, Ch. 1]. Notice that a group scheme embedding $\mu_p\hookrightarrow E^{\rm sm}$
gives data $C$ and $\gamma$; in fact, if $p$ is invertible on $S$, giving $C$ together with $\gamma$ as above
exactly amounts to giving a group scheme embedding $\mu_p\hookrightarrow E^{\rm sm}$.) We denote by $X_1 = X_1(p)$
the corresponding coarse moduli scheme over $\Spec(\Z)$.  The group $ (\Z/p\Z)^* $ acts on 
${\mathfrak M}_{\Gamma_1(p)}$ via 
\begin{equation}\label{action}
(a\, {\rm mod}\, p)\cdot (E, C,  \gamma)=(E, C,   \gamma\circ a^{-1}).
\end{equation}
(When $p$ is invertible, this action sends the corresponding $j: \mu_p\hookrightarrow E^{\rm sm}$ to 
the composition $j\circ (z\mapsto z^a): \mu_p\hookrightarrow E^{\rm sm}$.)
This produces a faithful $\Gamma= (\Z/p\Z)^*/\{\pm 1\}$-action on $X_1$.  When $H$ is a subgroup of $\Gamma  $ we let $X_H$ be the
quotient $X_1/H$, and define $X_H[1/p] = \Z[1/p] \otimes_{\Z} X_H$.  Set $X_0 = X_\Gamma$. 

The Tate curve $\bar {\mathcal G}_m/q^{\Z}$ over $\Spec(\Z[[q]])$ together
with the embedding $\mu_p\subset {\mathcal G}_m/q^{\Z}$ (see [DR, VII])
gives a morphism $\tau: \Spec(\Z[[q]])\to X_H$. We call the 
support of the corresponding section
$\Spec(\Z)\to \Spec(\Z[[q]])\to  X_H$ the $\infty$ cusp. 
Over $\C$, provided we trivialize $\mu_p(\C)$ via $\zeta_p=e^{2\pi i/p}$, this corresponds to the ``usual" $\infty$ cusp
and the parameter $q$ to $e^{2\pi iz}$ with $z$ in the upper half plane ${\mathfrak H}$.
The morphism $\tau$ identifies $\Spec(\Z[[q]])$
with the formal completion of
$X_H$ along $\infty$.

\begin{Remarknumb}\label{alternative}
{\rm a) Alternatively, we can define ${\mathfrak M}_{\Gamma_1(p)}$ by considering the non-proper moduli stack
of triples $(E,C,\gamma)$ as above with $E\to S$ an ``honest" (i.e smooth) elliptic curve
and then take its normalization over the affine $j$-line (compare to the approach in [DR], [KM]).

b) ${\mathfrak M}_{\Gamma_1(p)}$ is a variant of the following moduli stack which
is used in [CES]:
Let ${\mathfrak M}'_{\Gamma_1(p)}$ be the moduli stack
of triples $(E, \Sigma, \gamma')$ where $E\to S$ is a generalized elliptic curve, $\Sigma$
a locally free rank $p$ subgroup of the smooth locus $E^{\rm sm}$
which intersects every  irreducible component
of every geometric fiber of $E\to S$ and $\gamma': (\Z/p\Z)_S\to \Sigma$
a ``generator" of   $\Sigma$. If $(E,C,\gamma)$ is a point of  ${\mathfrak M}_{\Gamma_1(p)}$ 
with $E$ an honest elliptic curve, then 
the Weil pairing induces a canonical isomorphism $C^D\xrightarrow{\iota} E[p]/C$
($E[p]$ is the $p$-torsion subgroup scheme of $E$). The map
$(E,C,\gamma)\mapsto (E/C, E[p]/C,\gamma')$ with $\gamma':=\iota\,\circ\,\gamma$
induces an isomorphism of moduli stacks
on the complement of the cusps; this extends (by normalization) to an isomorphism
 ${\mathfrak M}_{\Gamma_1(p)}\xrightarrow{\sim} {\mathfrak M}'_{\Gamma_1(p)}$ which in turn induces
an isomorphism  between the corresponding coarse moduli 
schemes (compare [CES, p. 381]).  
}\end{Remarknumb}

Since we have assumed $p \equiv 1 $ mod $24$, the genus $g_0$ of $(X_0)_{\C}$ is $(p-13)/12$.
The work of Deligne and Rapoport [DR],
Katz and Mazur [KM] and Conrad, Edixhoven and Stein in [CES] implies
the following two results.

\begin{thm}\label{ramstruct} 

\begin{enumerate}
\item[a.] The scheme $X_H\to \Spec(\Z)$ is a flat projective curve, $X_H$ is normal Cohen-Macaulay and 
$X_H[1/p]\to \Spec(\Z[1/p])$ is smooth. The special fiber of  $X_H$ over $p$ has two irreducible
components $D^H_\infty$ and $D^H_0$ distinguished by the fact that $D^H_\infty$ intersects
the cuspidal section $\infty$; these have multiplicities $1$ and $(p-1)/(2\cdot\#H)$ respectively.
\item[b.]
The scheme $X_H$ has at most two non-regular points which are rational singularities
and lie on $ D^H_0-(D^H_0\cap D^H_\infty)$. 
Their exact number depends on $\# H\mod 6$: In particular, if $6$ divides  $\# H$ there are no 
such points and   $X_H$ is then regular. If $H=\{1\}$ there are two non-regular points on $X_1$.
There is a morphism $b: X'_1\to X_1$ which is a rational resolution of those two singular points
and a morphism $c:X'_1 \to \X_1$ which is a sequence
of blow-downs of exceptional curves such that $\X_1$ is regular and all the geometric
fibers of $\X_1 \to {\rm Spec}(\Z)$ are integral.  
Let $U=  X_1-D^{\{1\}}_0\subset X_1$. Then $U\to \Spec(\Z)$ is smooth, 
$b$ and $c$ are isomorphisms on $b^{-1}(U)$ and $\X_1 - c(b^{-1}(U))$ has dimension $0$.
\item[c.] The special fiber of $X_0=X_\Gamma$ over $p$
is reduced with simple normal crossings. Each of the two irreducible components 
$D_\infty=D^\Gamma_\infty$ and $D_0=D^\Gamma_0$ are isomorphic to ${\bf P}^1_{{\bf F}_p}$
and $D_0 \cdot D_\infty = g_0 + 1 = (p-1)/12$.
\end{enumerate}
\end{thm}

\begin{Proof}
Parts (a) and (b) follow from [CES, Th. 4.1.1, Th. 4.2.6]  (see
[CES, \S 5.3]; the reader should be aware that [CES] use the variant moduli problem
of Remark \ref{alternative} (b).) Part (c) follows from [DR, VI. 6.16]. \endproof
\end{Proof}

\begin{thm}\label{ramstruct2}
Assume that $6$ divides the order $\#H$. 
\begin{enumerate}
\item[a.]
 The morphism $\pi_H: X_H\to X_0$ is a tame $G = \Gamma/H$ cover of regular
projective curves and 
$\pi_H[1/p]:X_H[1/p]\to X_0[1/p]$ is a $G$-torsor. 
\item[b.] The morphism $\pi_H$ is totally ramified over the generic point 
of $D_0$, and unramified over the generic point of $D_\infty$.
The irreducible components $D^H_0$ and $D^H_\infty$ of $X_H\otimes_{\Z}{\fp}$ 
are the (reduced) inverse
images of $D_0$ and $D_\infty$ under $\pi_H$.   The character $\chi_{D^H_0}$
giving the action of $G$ on the  cotangent space of the codimension
$1$ generic point of  $D^H_0$ equals $\omega^{-2\cdot \#H}$, where 
$\omega: (\Z/p\Z)^*\to {\bf F}_p^*$ 
is the Teichm\"uller (identity) character. 
\end{enumerate}
\end{thm}

\begin{Proof}
The fact that $\pi_H[1/p]$ is a $G$-torsor follows from [Ma, II. 1];
the rest of part (a) then follows from Theorem \ref{ramstruct}. It remains to show part (b).
We can see from the references above that the geometric closed points of $D^{\{1\}}_0$, resp. $D^{\{1\}}_\infty$, 
which correspond to ordinary elliptic curves are given by triples
$(E, \Z/p\Z\subset E, 0)$, resp. $(E, \mu_p\subset E, \gamma)$ 
with $\gamma: \Z/p\Z\xrightarrow{\sim} \mu_p^D=\Z/p\Z$. Since by Theorem \ref{ramstruct} (a) $D^H_\infty$ has multiplicity $1$,
$\pi_H$ is unramified over the generic point of $D_\infty$. It remains to consider $D^H_0$;
our claim is a statement about
the completion of the local ring at the generic point of $D^H_0$.
It is enough to show this statement after a base change from $\Z$ to 
the Witt ring $W=W(\acfp)$ of an algebraic closure of $\fp$.

Let $E$ be an ordinary elliptic curve over $\acfp$.
The pair $(E, \Z/p\Z\subset E, 0)$ corresponds to a point $s$ of $D^{\{1\}}_0$. 
Let ${\mathcal R}_s$ be the formal deformation 
ring of the point $(E,  \Z/p\Z\subset E, 0)$ in the moduli stack ${\mathfrak M}_{\Gamma_1(p)}$.
Then ${\mathcal R}_s$ supports an action of $ \Delta\times (\Z/p\Z)^*$ where 
$\Delta={\rm Aut}(E)$. Let
$H'$ be the inverse image of $H$ under the surjection $(\Z/p \Z)^* \to \Gamma$,
and let $ s'$, $s''$ be the images of $s$ on $X_H$, resp. $X_0$.
The completions of the local rings are
\begin{equation}
\label{eq:quots}
\widehat {\O}_{X_H \otimes_{\Z} W,s'}\simeq ({\mathcal R}_s)^{\Delta\times H'}\quad {\rm and}\quad 
\widehat {\O}_{X_0 \otimes_{\Z} W,s''}\simeq ({\mathcal R}_s)^{\Delta\times (\Z/p)^*}
\end{equation}
as rings with $G=(\Z/p\Z)^*/H'$-action. 
 Using [CES, Theorem 3.3.3] and the proof of [CES, Theorem 4.1.1]
(or alternatively a direct calculation using 
the description of the scheme of generators of the multiplicative 
group scheme $\mu_p$ given in [KM, II])
we can deduce the following:
There is a $\Delta\times (\Z/p\Z)^*$-isomorphism
\begin{equation}
{\mathcal R}_s\simeq W[[v, u]]/(v^{p-1}-p) 
\end{equation}
where the action on the right hand side is as follows:
$\delta\in \Delta$ acts via 
$\delta\cdot v=\psi(\delta)v$, $\delta\cdot u=\psi^2(\delta)u$
with $\psi$ a faithful character of $\Delta$,  
while $a\in (\Z/p\Z)^*$ acts via 
$a\cdot v=\omega^{-1}(a)v$,
$a\cdot u=u$.
Assume now that $E$ is an ordinary elliptic curve over $\acfp$ with $j\neq 0, 1728$. 
Then $\Delta=\{\pm 1\}$. Part (b)  
for $D^H_0$ follows then directly from the above and (\ref{eq:quots}).
\endproof
\end{Proof}
\smallskip
 
  \begin{prop}\label{p1}
Let $H$ be an arbitrary subgroup of $\Gamma$. The group ${\rm H}^0(X_H, \O_{X_H})$ is isomorphic to $\Z$  with trivial $\Gamma/H$-action and ${\rm H}^1(X_H, \O_{X_H})$ is $\Z$-torsion free. 
\end{prop}
 
\begin{proof}
 We have ${\rm H}^0(X_H, \O_{X_H})\simeq \Z$
because $X_H\otimes_\Z\C$ is connected and 
$X_H\to \Spec(\Z)$ is projective and flat. 
The claim about ${\rm H}^1(X_H, \O_{X_H})$ follows
from Theorem \ref{ramstruct} (a) and [BLR, 9.7/1] (see 
Proposition \ref{lie}).
 \end{proof}
\smallskip

If $R$ is a subring of $\C$ we will denote by $S_2(\Gamma_1(p), R)$ the $R$-module
of cusp forms $F(z)=\sum_{n\geq 1}a_ne^{2\pi in    z}$ for the congruence subgroup
$\Gamma_1(p) \subset {\rm PSL}_2(\Z)$ 
whose Fourier coefficients $a_n$  belong to $R$.  (These are the Fourier coefficients ``at 
the cusp $\infty$".)   If $M$ is a finitely generated $\Z [\Gamma]$-module, let $M^\vee$ be the
$\Z [\Gamma]$ module $\Hom_{\Z}(M,\Z)$ with $\Gamma$ action $(af)(z) = f(a^{-1}z)$
for $a \in \Gamma$, $f \in M^\vee$ and $z \in M$.

\begin{prop}\label{p2}
There are $\Gamma$-equivariant isomorphisms 
\begin{equation*}
S_2(\Gamma_1(p), \Z)\simeq {\rm H}^0(X_1, \omega_{X_1/\Z})
\simeq {\rm H}^1(X_1, \O_{X_1})^\vee
\end{equation*}
where the $\Gamma$-action on $S_2(\Gamma_1(p), \Z)$ is via the diamond operators
and $\omega_{X_1/\Z}$ denotes the canonical (dualizing) sheaf of $X_1\to \Spec(\Z)$. 
\end{prop}

\begin{Proof}
The $\Gamma$-isomorphism 
\begin{equation*}
{\rm H}^0(X_1, \omega_{X_1/\Z})
\simeq {\rm H}^1(X_1, \O_{X_1})^\vee
\end{equation*}
follows from duality ([DR, I, (2.1.1)]).  Let $G(q)=\sum_{n\geq 1}a_nq^n\in S_2(\Gamma_1(p),\Z)$ with $q=e^{2\pi inz}$ and
consider $ G(q)dq/q$ as a regular differential over $\Spec(\Z[[q]])$.
A standard argument using the Kodaira-Spencer map shows that $G(q)dq/q$
extends to a regular differential over $X_1[1/p]$ (see for example [Ma, II \S 4]). 
This extension must also be regular in an open neighborhood of the section
at $\infty$.  Hence there is an open subset $U'$ of the set $U \subset X_1$
defined in Theorem \ref{ramstruct} (b) such that $G(q)dq/q$ is regular on $U'$
and $U - U'$ is a finite set of closed points.  
We obtain an injective $\Gamma$-equivariant homomorphism
\begin{equation}
\Phi: S_2(\Gamma_1(p),\Z)\to {\rm H}^0(U', \omega_{U'/\Z})={\rm H}^0(\X_1 , \omega_{\X_1/\Z})=
{\rm H}^0(X_1, \omega_{X_1/\Z})
\end{equation}
where the latter two equalities follow from the fact that $b:X'_1\to X_1$ and $c:X'_1 \to \X_1$ are rational morphisms
which are isomorphisms on $b^{-1}(U')$ and that $\X_1 - c(b^{-1}(U'))$ has codimension $2$ in $\X_1$.  
The surjectivity of $\Phi$ follows from pulling back elements of ${\rm H}^0(X_1,\omega_{X_1/\Z})$ via 
$\tau: \Spec(\Z[[q]])\to U$
and using the Kodaira-Spencer 
isomorphism.
\endproof
\end{Proof}
 
 \begin{prop}\label{p3}
We have $G=\Gamma/H$-isomorphisms
\begin{equation*}
S_2(\Gamma_H(p),\Z):=S_2(\Gamma_1(p),\Z)^H
\simeq {\rm H}^1(X_H, \O_{X_H})^\vee.
\end{equation*}
\end{prop}

\begin{Proof} Denote by $M_H$ the $H$-coinvariants of a $\Gamma$-module 
$M$. We then have a $\Gamma/H$-isomorphism $(M_H)^\vee\simeq (M^\vee)^H$.
In view of Proposition \ref{p2}, it will be enough to exhibit a 
$\Z[\Gamma/H]$-isomorphism  
\begin{equation}\label{want}
({\rm H}^1(X_1, \O_{X_1})_H)^\vee \simeq {\rm H}^1(X_H, \O_{X_H})^\vee.
\end{equation}
Let $\mu: X_1\to X_H$ be the quotient morphism.
Since taking $H$-coinvariants is a right exact functor and the 
coherent cohomology groups ${\rm H}^i(X_H, -)$ are trivial when $i\geq 2$
we obtain
\begin{equation}\label{cohcoinv}
{\rm H}^1(X_1, \O_{X_1})_H\simeq {\rm H}^1(X_H, \mu_*\O_{X_1})_H
\simeq {\rm H}^1(X_H, (\mu_*\O_{X_1})_H).
\end{equation}
The element $\Tr=\sum_{h\in H}h$
of the group ring $\Z[H]$ acts on $\mu_*\O_{X_1}$ to give a $\Gamma/H$-equivariant
morphism of  
$\O_{X_H}$-sheaves
\begin{equation}
(\mu_*\O_{X_1})_H\xrightarrow{\rm Tr} \O_{X_H}=(\mu_*\O_{X_1})^H.
\end{equation}
Since the stalks of $\mu_*\O_{X_1}$ are $H$-cohomologically trivial
at points where the ramification of $\mu:X_1\to X_H$ is tame (see [CEPT1]),
the cokernel $\mathcal C$ and kernel $\mathcal K$ of $\Tr$ are coherent and supported on the subset of 
$X_H$ over which the cover $\mu$ has wild ramification.  By [Ma,
II, \S 2] this is (at most) a finite set of points in characteristics
$2$ and $3$. Hence the cohomology groups of $\mathcal C$ and $\mathcal K$ vanish
in positive dimensions.  It follows that $\Tr$ induces a surjection
${\rm H}^1(X_H, (\mu_*\O_{X_1})_H )\to {\rm H}^1(X_H, \O_{X_H})$ with
finite kernel, so (\ref{want}) follows from  (\ref{cohcoinv}).
\end{Proof}

The following result is a corollary of a theorem of Rim [Ri].

\begin{prop}
\label{prop:rim} 
Let $\chi:  \Gamma \to \Z[\zeta_r]^*$ be a $1$-dimensional character 
of prime order $r\ge 5$ with kernel $H$. 
 Let $G = \Gamma/H$ and suppose $M$ is a finitely generated torsion-free $\Z[G]$-module.
Define $M^{\vee,\chi}$ to be the $\Z[\zeta_r]$-module $(M^\vee\otimes \Z[\zeta_r]\chi^{-1})^G$.   
\begin{enumerate}
\item[a.] There is a unique homomorphism $e'_\chi:\Gr_0(\Z[G]) \to \Cl(\Z[\zeta_r])$ 
such that for all $M$ as above, either $M^{\vee,\chi} = \{0\}$ and $e'_\chi([M]) = 0$ 
or $M^{\vee,\chi}$ is 
isomorphic to $\Z[\zeta_r]^s \bigoplus \mathfrak {U}$ for some integer $s \ge 0$ and a
$\Z[\zeta_r]$-ideal $\mathfrak{U}$ in the ideal class $e'_\chi([M])$.  
\item[b.] There is a unique isomorphism
$t_\chi:\Kr_0(\Z[G]) \to \Z \bigoplus \Cl(\Z[\zeta_r])$ such that
$t_\chi([P]) = ({\rm rank}_{\Z[G]}(P),e_\chi (\overline{[P]}))$ if $P$ is a projective $\Z[G]$-module,
where $\overline{[P]}$ is the image of $P$ in $\Cl(\Z[G])$ and $e_\chi:\Cl(\Z[G]) \to \Cl(\Z[\zeta_r])$ 
is the unique homomorphism such that $e_\chi(\overline{[P]}) = e'_\chi(f([P]))$ for all projective $P$, 
where $f:\Kr_0(\Z[G]) 
\to \Gr_0(\Z[G])$ is the forgetful homomorphism.
\end{enumerate}
\end{prop}

\medbreak
\noindent {\bf Proof of Theorem \ref{qua1}. }
\medbreak
With the notations of the Theorem, recall 
that  
$S_2(\Gamma_1(p),\Z[\zeta_r])_\chi$ is the 
$\Z[\zeta_r]$-submodule of $S_2(\Gamma_1(p),\Z[\zeta_r])$ consisting of cusp
forms of weight $2$ and of Nebentypus character $\chi$ 
whose Fourier coefficients at $\infty$ are in $\Z[\zeta_r]$.
Proposition \ref{p2} and its proof together with the fact that formation of 
the canonical sheaf commutes with the base change $\Z\to\Z[\zeta_r]$ implies that $S_2(\Gamma_1(p), \Z[\zeta_r])\simeq 
S_2(\Gamma_1(p), \Z)\otimes_\Z\Z[\zeta_r]$.  Propositions \ref{p2} and \ref{p3} 
now give
an isomorphism of (torsion free) $\Z[\zeta_r]$-modules
\begin{equation}
\label{eq:s2newisom}
S_2(\Gamma_1(p), \Z[\zeta_r])_\chi\simeq 
{\rm H}^1(X_H, \O_{X_H})^{\vee,\chi}. 
\end{equation}
(Here we are using the notation of Proposition \ref{prop:rim}.)  

The projective class $\chi^P(X_H,\O_{X_H}) \in \Kr_0(\Z[G])$ has the
property that
\begin{equation}
\label{eq:forget}
f(\chi^P(X_H,\O_{X_H}) ) = 
[{\rm H}^0(X_H, \O_{X_H})] - [{\rm H}^1(X_H, \O_{X_H})]
\end{equation} 
where  $f:\Kr_0(\Z[G]) \to \Gr_0(\Z[G])$ is the forgetful homomorphism.    If $P$
is a projective $\Z[G]$-module, then $\Q \otimes_{\Z} P$ is a free
$\Q[G]$-module, so ${\rm rank}_{\Z[\zeta_r]} (P^{\vee,\chi})  =
{\rm rank}_{\Z[G]}(P) = {\rm rank}_{\Z}(P)/r$.  Therefore (\ref{eq:forget}) implies
\begin{equation}
\label{eq:ranker}
{\rm rank}_{\Z[\zeta_r]} ({\rm H}^0(X_H, \O_{X_H})^{\vee,\chi} )
- {\rm rank}_{\Z[\zeta_r]} ({\rm H}^1(X_H, \O_{X_H})^{\vee,\chi}) = 
(1 - g(X_H))/r
\end{equation}
where $g(X_H)$ is the genus of $X_H$ and ${\rm H}^0(X_H,\O_{X_H}) \simeq \Z$
with trivial $G$-action by Proposition \ref{p1}.  Thus ${\rm H}^0(X_H, \O_{X_H})^{\vee,\chi}
= 0$ since $\chi$ is non-trivial.  Because the generic fiber of $X_H \to X_0$ is 
\'etale of degree $r$, we now conclude from (\ref{eq:s2newisom}) and 
the Hurwitz Theorem that
$$n(\chi) := {\rm rank}_{\Z[\zeta_r]}(S_2(\Gamma_1(p), \Z[\zeta_r])_\chi)  = 
(g(X_H) - 1)/r = g(X_0) - 1.$$
Since $p \equiv 1$ mod 24, we have $g(X_0) = (p-13)/12$, so 
$n(\chi) = (p - 25)/12$ as stated in Theorem \ref{qua1}.  

Recall that $\overline{\chi}^P(X_H,\O_{X_H})$ is the image of $\chi^P(X_H,\O_{X_H})$
in $\Cl(\Z[G])$.  We 
conclude from (\ref{eq:s2newisom}), (\ref{eq:forget}) and Propositions \ref{p1} 
and \ref{prop:rim} that there is
an isomorphism of $\Z[\zeta_r]$-modules
\begin{equation}
\label{eq:s2equal}
S_2(\Gamma_1(p), \Z[\zeta_r])_\chi \simeq \Z[\zeta_r]^{n(\chi) - 1} \bigoplus 
\mathfrak {U}
\end{equation}
where $\mathfrak{U}$ is a $\Z[\zeta_r]$-ideal having ideal class
$-e_\chi (\overline{\chi}^P(X_H,\O_{X_H}))$.

 We now compute the image of $\overline{\chi}^P(X_H,\O_{X_H})$ under the homomorphism
$\Theta = \Theta_3:\Cl(\Z[G]) \to C_{\Z}(G,3)$
by applying our main result, Theorem \ref{mainthm}, to the cover $\pi_H: X_H\to X_0$.
(Notice that $g(X_0)$ is odd; hence, we can use the version of the Theorem which does not include 
the factor of $2$.)
Since the index of $H$ in $\Gamma$ is the prime $r\geq 5$, the order of $H$ is divisible
by $6$. By Theorem \ref{ramstruct2} $\pi_H$ is ramified only at the fiber over $p$.
The field $\Q_p$
already contains a primitive $p-1$-st root of unity. Hence, we may take $R'_{(p)}=\Z_p$.
We find that
$\Theta(\bar\chi^P(X_H, \O_{X_H})))\in C_\Z(G;3)$ is given by the 
id\`ele $(b_v)_v\in {\bf A}^*_{f, \Q[G^3]}$
which is $1$ at all places $v\neq (p)$ and is such that
\begin{equation}
\label{eq:stu}
(\chi\otimes\phi \otimes\psi)(b_{(p)})=p^{- T((\chi-1)(\phi-1)(\psi-1))}
\end{equation}
with $T: {\rm Ch}(G)_{p}\to \Q$ the function  associated to the cover
$X_H\otimes_{\Z}{\Z_p}\to X_0\otimes_{\Z}\Z_p$ in Theorem \ref{mainthm}. 
For $a \in \Z/r\Z$ let $\{a\}$ be the unique integer in the range $0 \leq  \{a\} < r$ having residue class $a$.  By Theorem \ref{ramstruct}
(c) and Theorem \ref{ramstruct2}
the expression in (\ref{surf2}) becomes 
\begin{equation}\label{equ1}
T(\psi) = \frac{1-p}{12}\cdot \left ( \frac{g(\psi, D_{0})^2}{2} +\frac{g(\psi, D_{0})}{2}\right ) + g(\psi, D_{0}) =
\frac{1-p}{12} \left ( \frac{\{a\}_r^2}{2r^2} - \frac{\{a\}_r}{2r} \right ) - \frac{\{a\}_r}{r}.
\end{equation}
when $\psi = \chi_0^{-a}$, $\chi_0 =  \omega^{\frac{(p-1)}{r}}$ and $\omega:(\Z/p\Z)^* \to \Z_p^*$
is the Teichm\"uller character.

For $a \in \Z/r\Z$ define $\omega_r(a) = 0$ if $a = 0$, and otherwise let $\omega_r(a) \in \Z_r \subset \hat \Z$ be  the Teichm\"uller character associated to $r$.  Define
\begin{equation}
\label{eq:equ1}
T_1(\psi) =
\frac{1-p}{12} \left ( \frac{\omega_r(a)^2}{2r^2}\right ) \quad {\rm and} \quad 
T_2(\psi) = -\frac{\{a\}_r}{r}
\end{equation}
when $\psi = \chi_0^{-a}$ as above.  We extend $\psi \mapsto T_i(\psi)$ to a function on
the character ring ${\rm Ch}(G)_p$ by additivity.  
Since $p \equiv 1$ mod $24$ and $r|\frac{1-p}{24}$, we can define $\beta = (\beta_v)_v$
with $\beta_v \in \hat {\Z} \otimes_{\Z} \Z_v[G]^*$  by 
\begin{equation}\label{betaidele}
\psi(\beta_v)=\begin{cases}
1, &\text{if $v\neq (p)$;}\\
p^{ -T(\psi) + T_1(\psi) + T_2(\psi)}, 
&\text{if $v=(p)$.\ \ }
\end{cases}
\end{equation}
Since $\Cl(\Z[G])$ is a torsion group, $\beta$ defines a unique class $[\beta]$ in
$\Cl(\Z[G])$.  

We now show 
\begin{equation}
\label{eq:almostP}
\overline{\chi}^P(X_H,\O_{X_H}) = [\beta].
\end{equation}

Define $D = [\beta]-\overline{\chi}^P(X_H,\O_{X_H})$, and let
$R = \hat{\Z}$ if $i = 1$ and $R = \Z$ if $i = 2$.  From
(\ref{eq:equ1}) one has $rT_i(\psi) \in R$ and $rT_i(\psi)  \equiv a$ mod $rR$.
It follows that for all triples $(\chi,\phi,\psi)$ elements of $Ch(G)_p$, $T_i(\chi\phi\psi-\chi\phi-\phi\psi-\chi\psi+\chi+\phi+\psi-1)$ lies in $R$.  
Hence there are  
elements $a_i = (a_{i,v})_v \in \prod_{v\ {\rm finite}}
(R \otimes_{\Z} \Q_v[G^3]^*)$ for which
\begin{equation}\label{modformprime}
(\chi\otimes\phi\otimes\psi)(a_{i,v})=\begin{cases}
1, &\text{if $v\neq (p)$;}\\
p^{ T_i(\chi\phi\psi-\chi\phi-\phi\psi-\chi\psi+\chi+\phi+\psi-1)}, 
&\text{if $v=(p)$.\ \ }
\end{cases}
\end{equation}
 We now conclude from 
 (\ref{eq:stu}), (\ref{equ1}) and (\ref{betaidele}) 
that 
\begin{equation}
\label{eq:thetadeq}
1 \otimes \Theta(D) = c_1 + c_2
\end{equation}
in $\hat {\Z} \otimes_{\Z} C_{\Z}(G;3)$, 
where $c_i$ is the class associated to the element $a_i$.

Let us first show
\begin{equation}
\label{eq:c2iszero}
c_2 = 0.
\end{equation}
For this it will suffice to show that there is a cubic element
$\lambda \in \Q[G^3]^*$ such that $\lambda a_2$ is a unit
id\`ele of $\Q[G^3]$.  Fix a primitive $p$-th root of unity
$\zeta_p \in \overline{\Q}^*$, and let 
$$\tau(\psi) = \sum_{j \in (\Z/p)^*} \psi(j) \zeta_p^j$$
be the usual Gauss sum associated to $\psi$.  Let
$\tau$ be the unique extension of the map $\psi \mapsto \tau(\psi)$
to a homomorphism from $R_G$ to ${\overline {\Q}}^*$.   We let
$\tau^{(3)}$ be the element of $\mathrm{Hom}(R_{G^3},{\overline{\Q}}^*)$
which sends $(\chi,\phi,\psi)$ to 
$$\tau(\chi\phi\psi-\chi\phi-\phi\psi-\chi\psi+\chi+\phi+\psi-1).$$
From the behavior of Gauss sums under automorphisms of $\overline{\Q}$,
and the factorization of the ideals they generate (c.f. [La1, \S IV.3]), it follows
that 
$\tau^{(3)}$ is $\mathrm{Gal}(\overline{\Q}/\Q)$-equivariant,
and corresponds to an element $\lambda \in \Q[G^3]^*$
of the required kind.  This shows (\ref{eq:c2iszero}).

Turning now to $c_1$, let $\sigma(s)$ be the automorphism of $G$
which sends $g \in G$ to $g^s$ for $s \in(\Z/r\Z)^*$.  By (\ref{idel2}) the action of ${\rm Aut}(G)$
on $\Z_p[G]^*$ corresponds to the  action of
${\rm Aut}(G)$ on $f \in \Hom({\rm Ch}(G)_p,\Q_p^*) $
defined by $(\sigma(s)( f))(\chi) = f(\sigma(s)^{-1}(\chi))  = f(\chi^s)$
for $\chi \in {\rm Ch}(G)_p$.  
From the definition of the $T_1$ in (\ref{eq:equ1}) and the multiplicativity of
the Teichm\"uller character we have
$(\sigma(s)T')(\psi) = \omega_r(s)^2 T'(\psi)$.  
It follows that the element $\alpha(s) = \sigma(s) - \omega_r(s)^2$ 
 sends $a_1$ to the identity function, so 
 \begin{equation}
 \label{eq:handlec1}
 \alpha(s) c_1 = 0
 \end{equation}
 
 Because
$\Theta$ is ${\rm Aut}(G)$-equivariant, we can now conclude from 
(\ref{eq:thetadeq}), (\ref{eq:c2iszero}) and (\ref{eq:handlec1}) that
\begin{equation}
\label{eq:killer}
1 \otimes \left ( \Theta(\alpha(s) \cdot D) \right ) = 0\quad {\rm in}\quad 
\hat{\Z} \otimes_{\Z} C_{\Z}(G;3).
\end{equation}
If $A \to B$ is an injection of abelian groups, and $A$ is finite,
then $A = \hat{\Z}\otimes_{\Z} A \to \hat{\Z}\otimes_{\Z} B$ is injective,
as one sees by reducing to the case in which $A\to B$ is the inclusion
$n^{-1}\Z/\Z \to \Q/\Z$ for some $n \ge 1$.  So (\ref{eq:killer})
and the injectivity of $\Theta$ implies 
\begin{equation}
\label{eq:killer2}
\alpha(s) \cdot D = 0\quad {\rm in}\quad 
\Cl(\Z[G]).
\end{equation}
Similarly, since $r^2 T_1(\psi) $ is in ${\Z}_r \subset \hat {\Z}$
and $\Cl(\Z[G])$ is a torsion group, we see from the injectivity of
$\Theta$ that $D$ is in the $r$-Sylow subgroup of $\Cl(\Z [G])$.  

We now use the fact that $\Cl(\Z[G])$ is isomorphic to $\Cl(\Z[\zeta_r])$. 
 Define
$C_{j}$ to be the group of classes $c$ in the $r$-Sylow
subgroup of $\Cl(\Z[\zeta_{r}])$ for which
$(\sigma(s) - \omega_r(s)^j)(c) = 0$
for all $s \in(\Z/r)^*$.  We have shown $a_1$ corresponds to a class in $ C_{2}$.   
By the Spiegelungsatz (c.f. [Wa, Theorem 10.9]),
$C_{2} = 0$ if
$C_{-1} = 0$.  Herbrand's Theorem ([Wa, Thm. 6.17])
shows that if $C_{-1} \ne 0$, then the Bernoulli number
$B_{r -(r-2)} = B_2$ is congruent to $0$ mod  $r{\bf Z}_r$.  This is impossible since $B_2 = \frac{1}{6}$ and $r \ge 5$, so we have shown (\ref{eq:almostP}).

In view of (\ref{eq:s2equal}) and (\ref{eq:almostP}), the proof of Theorem \ref{qua1} is reduced to showing
\begin{equation}
\label{eq:zipit}
-e_\chi ([\beta]) = \theta_2\cdot [\PP_\chi]
\end{equation}
when $\theta_2$ and $\PP_\chi$ are as defined in Theorem \ref{qua1}.
By applying a suitable element of $\Delta = {\rm Gal}(\Q(\zeta_r)/\Q)$ to both
sides, we can reduce to the case in which $\chi = \chi_0 = \omega^{\frac{(p-1)}{r}}$.
From the definition of $e_\chi$ in Proposition \ref{prop:rim}, $\beta$ 
in (\ref{betaidele}) and of $\theta_2$ and $\PP_\chi$ we see that
$$-e_{\chi_{0}}([\beta]) = \theta_2 \cdot [\PP_{\chi_0}]  -  \frac{ (p-1)}{24r} \theta_1 \cdot [\PP_{\chi_0}]$$
$\theta_1 =  \sum_{a  \in (\Z/r)^*}   \{a\}  \sigma^{-1}_a \in \Z[\Delta]$.
By Stickelberger's Theorem, $\theta_1$ annihilates $\Cl(\Z[\zeta_r])$, so 
the proof is  complete.

\begin{Remarknumb}\label{generalp}{\rm Our techniques can be used to obtain a
result for cusp forms of (general) square free level $N$. In the general case, the moduli 
theoretic model for $X_0(N)$ is not regular over $\Z$. Hence, as in Remark \ref{abhya} (b),
we need to deal with  the cover
obtained by the
normalization of the modular curve $(X_1(N))_\Q$ over a regular resolution 
of $X_0(N)$. 
The details, which are somewhat involved, will appear in future work.
}
\end{Remarknumb}
 \bigskip

 \section{An equivariant Birch and Swinnerton-Dyer relation} \label{BSD}
 
Suppose that $G$ is a finite abelian group and let $V\to W$ be a $G$-cover 
of smooth projective geometrically connected curves of positive genus over $\Q$. We will denote by
$A={\rm Jac}(V)$ the Jacobian variety of $V$ (an abelian variety over $\Q$); the group $G$ acts
on $A$ by (Picard) functoriality. Let $A(\Q)$ and $\Sha (A)$ denote the Mordell-Weil group and the Tate-Shafarevich group of $A$
over $\Q$. 
%Both $A(\Q)$ and $\Sha (A)$ are naturally $\Z[G]$-modules. 
%By the Mordell-Weil theorem $A(\Q)$ is finitely 
%generated. In what follows, we will always assume that $\Sha(A)$ is finite.  
Let ${\rm H}_1(A(\C),\Z)^+$ be the subgroup of the (singular) homology group ${\rm H}_1(A(\C),\Z)$ which is fixed by the action of complex conjugation on $A(\C)$.
Define $\A$ to be the N\'eron model of $A$ over $\Spec(\Z)$; this 
supports a natural $G$-action that extends the action of $G$ on $A={\A}\times_\Z\Q$.
Consider    the tangent space ${\rm Lie}(\A)$ of $\A$ at the origin
section $e:\Spec(\Z)\to \A$. ${\rm Lie}(\A)$ is  a finitely generated $\Z[G]$-module which is
  free as a $\Z$-module.
For a prime number $l$ we denote by $ \Phi_l(A)$ the group
of ${\bf F}_l$-points of the (finite) group scheme of connected 
components of the fiber of $\A$ over ${\bf F}_l$.
We also denote by $\Phi_\infty(A)$ the group of connected components 
$A({\bf R})/A({\bf R})^0$ of the Lie group $A({\bf R})$.  
Suppose $R$ is a Dedekind ring.  Denote by $x\mapsto x^\vee$ the involution on ${\rm G}_0(R[G])$ defined in the
following way.
If $M$ is a finitely generated $R[G]$-module, there is a short exact sequence of 
$R[G]$-modules
\begin{equation}
0\to L_1\to L_0\to M\to 0
\end{equation}
in which $L_0$ and $L_1$ are $R$-torsion free. We set $[M]^\vee=[\Hom_R(L_0, R)]-[\Hom_R(L_1,R)]$.
(Here the $G$-structure on $\Hom_R(L, R)$ is given as usual by $(g\cdot f)(l)=f(g^{-1}l)$.)

\begin{conjecture}\label{BSDcon}
Assume that the Tate-Shafarevich group $\Sha(A)$ is finite. Then the identity
\begin{equation}\label{BSDeq}
[{\rm Lie}(\A)]-[{\rm H}_1(A(\C),\Z)^+]+[\Phi_\infty(A)]+\sum_l[\Phi_l(A)]+[A(\Q)]^\vee-[A(\Q)]+[\Sha(A)]=0
\end{equation}
is true in the Grothendieck group  ${\rm G}_0(\Z[G])$ of finitely generated $\Z[G]$-modules.
\end{conjecture}

\begin{Remarknumb}\label{bigrem}
{\rm  We call this a Birch and Swinnerton-Dyer relation because it should
follow from an appropriate ``refined" version of the classical 
Birch and Swinnerton-Dyer conjecture for abelian varieties. A trivialization of
the virtual module corresponding to the left hand side of (\ref{BSDeq}) should
arise from the leading
terms in the Taylor expansion at $s = 1$ of L-series associated to the $G$ action
on $A$, the height pairing on $A(\Q)$ and the period map of $A$. 
A prototype of  such an extension of the Birch and Swinnerton-Dyer conjecture 
can be found in the works of Gross (see [G], [BG]). 
Following the seminal work of Bloch and Kato, very general equivariant Tamagawa number conjectures 
for motives with endomorphisms have been developed by Fontaine and Perrin-Riou
and by Burns and Flach.  These include, in particular, refined 
versions of the 
Birch and Swinnerton-Dyer conjecture; see [BF] or the survey [F]. 
Conjecture \ref{BSDcon} should   be 
a consequence of more precise conjectural formulas of this kind,
though the details have not appeared in the literature;
we are grateful to Matthias Flach for a private communication on this matter.
}
\end{Remarknumb}

\begin{prop}
\label{lie}  
Let $X$ be a normal flat projective curve over $\Z$ on which $G$
acts and which has rational singularities.  Suppose that the general fiber $X_{\Q}$ is $G$-isomorphic to $V$,
and that the greatest common divisor of the multiplicities 
of the irreducible components of each geometric fiber of $X\to \Spec(\Z)$ is $1$.  Then ${\rm H}^1(X,\O_X)$ is 
$G$-isomorphic to ${\rm Lie}(\A)$ (and hence is $\Z$-torsion free), and ${\rm H}^0(X,\O_X) \cong \Z$ with trivial
$G$-action.
\end{prop}

\begin{Proof} Since $V$ is smooth and geometrically connected and $X\to \Spec(\Z)$ is 
projective and flat
we have
${\rm H}^0(X,\O_X) \cong \Z$ with trivial
$G$-action. The claim that ${\rm H}^1(X,\O_X)$ is 
$G$-isomorphic to ${\rm Lie}(\A)$ follows, under our assumptions,  from 
[BLR, 9.7/1] and its proof. (For a more general 
treatment of the relation between ${\rm Lie}(\A)$ and ${\rm H}^1(X,\O_X)$ see [LLR].) \endproof 
\end{Proof}

\begin{cor}\label{BSDcon2}
Assume that the Tate-Shafarevich group $\Sha(A)$ is finite. With the hypotheses of Proposition \ref{lie},
Conjecture \ref{BSDcon} predicts
\begin{eqnarray}\label{BSDeq2}
\ \ \ \ \ [\Z]-f([\chi^P(X,\O_{X})])  &=& [{\rm H}_1(A(\C),\Z)^+] -[\Phi_\infty(A)]-\sum_l[\Phi_l(A)] \\
&{\ }&- \ ([A(\Q)]^\vee)+[A(\Q)]-[\Sha(A)]\nonumber
\end{eqnarray}
in ${\rm G}_0(\Z[G])$, where $f:{\rm K}_0(\Z[G]) \to {\rm G}_0(\Z[G])$
is the forgetful homomorphism. 
\end{cor}

\subsection{Modular curves.}

In this section we specialize Corollary \ref{BSDcon2} to the modular curve
case considered in  \S \ref{s:modforms}.  Let $p$, $r \ge 5$ be primes such that
$p \equiv 1$ mod $24r$, and let $\Gamma = (\Z/p\Z)^*/\{\pm 1\}$
be the group of diamond operators acting on the model $X_1(p)$ of
the modular curve $X_1(p)_\Q$ described in \S \ref{s:modforms}.
By Theorem \ref{ramstruct} (a) $X=X_1(p)$ satisfies the hypotheses of
Proposition \ref{lie}.
The next statement follows now from Theorem \ref{ramstruct2}.

\begin{prop}
\label{prop:shimura}
Let $Sh$ be the subgroup of order $6$ in $\Gamma$.  For all subgroups $H$ containing
$Sh$, the cover $X_H=X_1(p)/H \to X_0 = X_1(p)/\Gamma$ is a tame $G= \Gamma/H$
cover of curves over $\Z$ which satisfies the hypotheses of  Theorem \ref{mainthm}.
Theorem \ref{mainthm} leads to an exact expression for the term $\chi^P(X_H,\O_{X_H}) \in \Kr_0(\Z[G])$
appearing in the prediction (\ref{BSDeq2}) when $X = X_H$.\endproof
\end{prop}

\begin{Example}{\rm Suppose $H = Sh$.  Then $X_{Sh} = X_1(p)/Sh$ is a model of the Shimura
 cover of $X_0(p)_\Q$.  This is the largest quotient of $X_1(p)_\Q$
 which is unramified over $X_0(p)_\Q$.  The formula
 for $\chi^P(X_{Sh}, \O_{X_{Sh}})$ resulting from Theorem \ref{mainthm} requires 
 using the injection $\Theta: \Cl(\Z[G]) \to C_\Z(G;3)$.  }
 \end{Example}
 
 In order  to use the more explicit formula of Theorem \ref{qua1}, we will
 now let $H$ be the subgroup of (prime) index $r$ in $\Gamma$. Our goal
 is to show Theorem \ref{qua2} and Corollary \ref{corBSD} of \S \ref{s:intro}.
 
 Fix a character
 $\chi:G \to \C^*$ of order $r$.
Tensoring $\Z[G]$-modules with the 
ring homomorphism $\Z[G]\to \Z[\zeta_r,\frac{1}{2r}]$  induced by $\chi$
gives a Steinitz class homomorphism 
$$
s_\chi:\Gr_0(\Z[G]) \to \Gr_0(\Z[\zeta_r,\frac{1}{2r}])/\{{\rm free \ modules}\} = 
\Cl(\Z[\zeta_r,\frac{1}{2r}])= \Cl(\Z[\zeta_r,\frac{1}{2}]]).$$

 \begin{lemma}
 \label{lem:easy}
 Suppose that $V = \Q \otimes_\Z X_H$ in Conjecture \ref{BSDcon}.
Then $s_\chi$ sends the classes $[{\rm H}_1(A(\C),\Z)^+]$, $[\Phi_\infty(A)]$, and $[\Phi_l(A)]$
for all primes $l$, to $0$
 in $ \Cl(\Z[\zeta_r,\frac{1}{2}])$.
\end{lemma}

\begin{Proof}  By [CES, Theorem 1.1.3, Corollary 1.1.5], $\Phi_p(A)$ is
a finite cyclic group with trivial action by $G$, so $s_\chi(\Phi_p(A)) = 0$.  Since $X_H$ has good
reduction outside of $p$, $\Phi_l(A) = \{0\}$ for $l \ne p$.  By [CES, Proposition 6.1.12],
$\Phi_\infty$ is a finite two-group, so $s_\chi([\Phi_\infty]) = 0$.
Finally, ${\rm H}_1(A(\C),\Z)$ is isomorphic to ${\rm H}_1(V(\C),\Z)$ as a module
for $G\times T$ where $T = \{e,c\}$ is the group generated by the action
of complex conjugation $c$ on $A(\C)$ and $V(\C)$.  Since the action of $G$
on $V(\C)$ is free, we can find a finite triangulation
of $V(\C)$ which is stable under the action of $G \times T$, and such
that the stabilizer of any element of the triangulation lies in $\{e\} \times T$.
Computing ${\rm H}_*(V(\C),\Z)$ using this triangulation, we see that
the Euler characteristic
\begin{equation}
\label{eq:betti}
[{\rm H}_0(V(\C),\Z))] - [{\rm H}_1(V(\C),\Z)] + [{\rm H}_2(V(\C),\Z)]
\end{equation}
in $\Gr_0(\Z[G \times T])$ is an integral combination of classes of permutation modules
of the form $\Z[(G \times T)/J]$ where $J \subset \{e\} \times T$.
The morphism $M \to M^+_2 := \Z[\frac{1}{2}]\otimes_\Z M^T$ is an exact
functor from the category of $\Z[G \times T]$-modules to the category
of $\Z[\frac{1}{2}][G]$-modules.  Applying this functor to (\ref{eq:betti})
and using the fact that ${\rm H}_0(V(\C),\Z)$ and ${\rm H}_2(V(\C),\Z)$
have trivial $G$-action leads to the conclusion that $s_\chi([{\rm H}_1(A(\C),\Z)^+])$
is a sum of classes of the form $s_\chi(M^+_2)$ where $M$ is a permutation module
of the form $\Z[(G\times T)/J]$ for some $J \subset \{e\} \times T$. One
checks readily that all such $s_\chi(M^+_2)$ are trivial,
which completes the proof.
\endproof
\end{Proof}

\medbreak
\noindent {\bf Proof of Theorem \ref{qua2} of \S \ref{s:intro}.}

With the notations of Lemma \ref{lem:easy} and Theorem \ref{qua1}, we are to show that Conjecture \ref{BSDcon} implies
\begin{equation}
%\label{bsd3}
 \overline{\theta_2[\PP_\chi]} = s_\chi( [\Sha(A)])-\overline{s_\chi( [A(\Q)])}-
s_\chi( [A(\Q)]) \label{firstline} 
\end{equation}
in $\Cl(\Z[\zeta_r, 1/2])$, where $A = \mathrm{Jac}(X_{H,\Q})$.

By Corollary \ref{BSDcon2}, Conjecture \ref{BSDcon} implies the equality (\ref{BSDeq2}). 
The equality (\ref{BSDeq2}) gives (\ref{firstline})
in view of (\ref{eq:almostP}), (\ref{eq:zipit}), Proposition \ref{prop:rim}, Lemma \ref{lem:easy}, $s_\chi([\Z]) = 0$  and the fact that 
\begin{equation}
\label{eq:dualtime}
s_\chi(M^\vee) = - \overline {s_\chi(M)}
\end{equation}
if $M^\vee = {\rm Hom}_\Z(M,\Z)$ and $\overline {\mathcal{D}}$ is the complex conjugate
of an ideal class $\mathcal{D}$.\endproof
\medbreak

Now let $C(p)$ be the subgroup of $J_1(p)(\Q)$ which is generated by 
the differences of $\Q$-rational cusps on $X_1(p)_\Q$.
(With our definition of $X_1(p)_\Q$, these are the cusps that lie above 
the cusp $\infty$ of $X_0(p)_\Q$).  By the Manin-Drinfeld Theorem,
$C(p) \subset J_1(p)(\Q)_{\mathrm{tor}}$.
As in the introduction, denote by $J'_H(\Q)$  the image of 
$A(\Q) $ in $J_1(p)(\Q)/C(p)$ under the 
homomorphism induced by the pullback map $\pi^*:A(\Q) \to J_1(p)(\Q)$.    

\begin{prop}\label{newpropschi}
Under the above assumptions, 
$$
s_\chi([J'_H(\Q)]) = s_\chi([A(\Q)])
$$
in $\Cl(\Z[\zeta_r, 1/2])$.
\end{prop}
\begin{Proof} 
The main ingredient is the following:

\begin{lemma}
\label{lem:cuspdivcl} 
The quotient morphism $\pi:X_1(p)_\Q \to V = (X_H)_\Q$ induces homomorphisms $\pi^*:A(\Q) \to  
J_1(p)(\Q)$ and $\pi_*:J_1(p)(\Q) \to A(\Q)$.  Let $T_H$ be
the trace element $\sum_{h \in H} h$ of the group ring $\Z[H]$.
\begin{enumerate}
\item[a.]
The composition $\pi^* \circ \pi_*$ (resp. $\pi_* \circ \pi^*$) is multiplication
by $T_H$ (resp. $\# H$) on $J_1(p)(\Q)$ (resp. $A(\Q)$).  
\item[b.] Both $\mathrm{kernel}(\pi^*)$ and the 
Tate cohomology group $\hat{\Hr}^0(H,C(p))$ are finite cyclic groups
with trivial action by $G$. 
\item[c.] The $G$-module $T_H(C(p))$ has trivial class in $\Gr_0(\Z[G])$. 
\end{enumerate}
\end{lemma}

\begin{Proof} Statement (a) is true because $\pi^*$ (resp. $\pi_*$)
is induced by the pullback (resp. pushdown) of divisors via
the $H$-cover $\pi$.  The morphism $\pi$ is a composition 
 $ X_1(p)_\Q\to (X_{Sh})_\Q\to V$ where the first, resp. second morphism, is
totally ramified, resp. unramified. It follows that the 
kernel of $\pi^*_{\overline \Q}$ is contained in the
kernel of $(\pi^*_{Sh})_{\overline \Q}: A(\bar\Q)\to {\rm Jac}(X_{Sh})(\bar\Q)$.
By Kummer theory, the field 
$\overline{\Q}(X_{Sh})$ is isomorphic to $\overline{\Q}(V)(f^{1/\delta})$
for some $f \in\overline{\Q}(V)^*$, where $\delta$ is the degree of $\pi_{Sh}$.  There is a divisor $c$ on $\overline{\Q} \otimes_\Q V$ such that $\delta c = \mathrm{div}(f)$, and the class $[c]$ of $c$ in $A(\overline{\Q})$
generates the kernel of $(\pi^*_{Sh})_{\overline \Q}$. Since $\Gamma$ is abelian, the  group $G$ acts
trivially on $[c]$. Therefore,  $\mathrm{kernel}(\pi^*)\subset \mathrm{kernel}(\pi^*_{\overline \Q})
\subset \mathrm{kernel}((\pi^*_{Sh})_{\overline \Q})$
has the properties stated in part (b).

Concerning $\hat{\Hr}^0(H,C(p))$ and $T_H (C(p))$, let 
$I_\Gamma$ be the augmentation ideal of $\Z[\Gamma]$.
Define $\mathcal{J}$ to be the kernel of surjection
$\Z[\Gamma] \to \mu_p^{\otimes 2}$ induced by the square of the Teichm\"uller
character $\omega_p:\Gamma \to \mu_p$.  In [KL, Thm. 3.4] it
is shown that there is an exact sequence of $\Z[\Gamma]$-modules
\begin{equation}
\label{eq:aargh}
0  \rightarrow I_\Gamma \cap \mathcal{J} \xrightarrow{t} I_\Gamma \rightarrow C(p) \rightarrow  0
\end{equation}
in which $t$ is multiplication by a certain second Stickleberger
element in $\Q[\Gamma]$ whose definition we will not require.  
From the definition of $\mathcal{J}$ and $p > 2$ we see there is an exact sequence
\begin{equation}
\label{eq:fiddle}
0 \to I_\Gamma \cap \mathcal{J} \xrightarrow{i} I_\Gamma \to \mu_p^{\otimes 2} \to 0.
\end{equation}
in which $i$ is inclusion.
Suppose that 
\begin{equation}
\label{eq:generalh}
0  \rightarrow I_\Gamma \cap \mathcal{J} \xrightarrow{h} I_\Gamma \rightarrow M \rightarrow  0
\end{equation}
is any exact sequence of $\Z[\Gamma]$-modules.  Taking Tate cohomology with 
respect to $H$ gives an exact sequence  
\begin{equation}
\label{eq:yupM}
\hat{\Hr}^0(H,I_\Gamma) \to \hat{\Hr}^0(H,M) \to {\Hr}^1(H,I_\Gamma \cap \mathcal{J}) \to 
{\Hr}^1(H,I_\Gamma)
\end{equation}
From the cohomology of the exact sequence $0 \to I_\Gamma \to \Z[\Gamma] \to \Z \to 0$ we
see that $\hat{\Hr}^0(H,I_\Gamma) \cong \hat{\Hr}^{-1}(H,\Z) = (0)$ and ${\Hr}^1(H,I_\Gamma) \cong \hat{\Hr}^0(H,\Z)$  is a finite
cyclic group with trivial $G$-action.  Hence  all submodules of 
${\Hr}^1(H,I_\Gamma)$ are finite and cyclic with trivial $G$-action.
When we let $(M,h) = (\mu_p^{\otimes 2},i)$ in (\ref{eq:yupM}) we get $\hat{\Hr}^0(H,M) = 0$
since $M^H = (\mu_p^{\otimes 2})^H = (0)$ because $\omega_p^{ 2}$ is non-trivial
on $H$. Therefore (\ref{eq:yupM}) shows ${\Hr}^1(H,I_\Gamma \cap \mathcal{J})$
is finite and  cyclic with trivial $G$-action.

We now apply
(\ref{eq:yupM}) with $(M,h) = (C(p),t)$. Since $\hat{\Hr}^0(H,I_\Gamma) = 0$, $\hat{\Hr}^0(H,C(p))$
is isomorphic to a  submodule of ${\Hr}^1(H,I_\Gamma \cap \mathcal{J})$. Hence $\hat{\Hr}^0(H,C(p))$ is finite and cyclic with trivial $G$-action and
$
[\hat{\Hr}^0(H,C(p))] = 0 \quad \mathrm{in}\quad G_0(\Z[G])
$.
This completes the proof of part (b) of Lemma \ref{lem:cuspdivcl}, and it
reduces the proof of part (c) to showing that
\begin{equation}
\label{eq:moretiring}
[C(p)^H] = 0 \quad \mathrm{in}\quad G_0(\Z[G])
\end{equation}
because $\hat{\Hr}^0(H,C(p)) = C(p)^H/T_H(C(p))$.

Taking the $H$-cohomology of (\ref{eq:generalh}) gives an exact sequence
\begin{equation}
\label{eq:reallydumb}
0  \rightarrow (I_\Gamma \cap \mathcal{J})^H \xrightarrow{h} (I_\Gamma)^H \rightarrow M^H \rightarrow  {\Hr}^1(H,I_\Gamma \cap \mathcal{J})
\end{equation}
We have already shown ${\Hr}^1(H,I_\Gamma \cap \mathcal{J})$ is finite and cyclic
with trivial $G$-action, so every submodule of this module has 
 trivial class in $G_0(\Z[G])$.  Therefore (\ref{eq:reallydumb}) implies
\begin{equation}
\label{eq:oyveybigtime}
[M^H] = [(I_\Gamma)^H] - [(I_\Gamma \cap \mathcal{J})^H] \quad \mathrm{in}\quad G_0(\Z[G]).
\end{equation}
Applying this with the pairs $(M,h)$ given by $(C(p),t)$ and $(\mu_p^{\otimes 2},i)$ gives
\begin{equation}
\label{eq:equalatlast}
[C(p)^H] = [(\mu_p^{\otimes 2})^H] \quad \mathrm{in}\quad G_0(\Z[G]).
\end{equation}
Hence (\ref{eq:moretiring}) follows from (\ref{eq:equalatlast}) and the fact previously
observed that $(\mu_p^{\otimes 2})^H$ is trivial.
\endproof
\end{Proof}
\smallskip

Now let us complete the proof of Proposition \ref{newpropschi}.
When $T_H$ is the trace element of $\Z[H]$, we have 
$T_H (C(p)) = \pi^* \pi_* C(p) \subset \pi^* A(\Q)$ by
Lemma \ref{lem:cuspdivcl}(a).  Hence
$W: = (C(p) \cap \pi^* A(\Q))/T_H(C(p))$ is contained in
$\hat{\Hr}^0(H,C(p)) = C(p)^H/T_H(C(p))$, so $W$ is cyclic
with trivial $G$ action by Lemma \ref{lem:cuspdivcl}(b).  
In $\Gr_0(\Z[G])$ we have
$$[J'_H(\Q)] + [W] + [T_H(C(p))] = [J'_H(\Q)] + [C(p) \cap \pi^* A(\Q)] = A(\Q) - [\mathrm{kernel}(\pi^*)].$$
Since $W$ and $\mathrm{kernel}(\pi^*)$  are finite cyclic with trivial $G$-action,
$s_\chi([W]) = s_\chi([\mathrm{kernel}(\pi^*)]) = 0$. We also have $[T_H(C(p))] = 0$
by Lemma \ref{lem:cuspdivcl}(c).  So $s_\chi([J'_H(\Q)]) = s_\chi([A(\Q)])$. \endproof
\end{Proof}

\begin{Remarknumb}\label{AV}
 {\rm We can also consider the following variant of the conjectural identity 
 in Theorem \ref{qua2}. The $\Gamma=(\Z/p\Z)^*/\{\pm 1\}$-action 
on $A={\rm Jac}(X_1(p)_\Q)$ induces
an embedding of $\Q[\Gamma]\simeq\oplus_{r}\Q(\zeta_r)$ 
in the endomorphism algebra ${\rm End}^0(A)={\rm End}(A)\otimes_\Z\Q$. The summands 
here correspond to 
characters $\chi: \Gamma\to \Z[\zeta_r]^*$ of order $r$. 
The  component of $\Q[G]$ that corresponds to $\chi$ determines an isogeny class $\{A_\chi\}$
of abelian varieties over $\Q$ with the following property: For each prime $l$ 
the (well-defined) ${\rm Gal}(\bar\Q/\Q)$-representation 
${\rm H}^1_{et}(A_{\chi}\otimes_\Q\bar\Q, \Q_l)$ is the 
$\chi$-component $V^\chi_l$ of $V_l={\rm H}^1_{et}(A\otimes_\Q{\bar \Q}, \Z_l)\otimes_{\Z_l}\Q_l$.
The field $\Q(\zeta_r)$ embeds into the endomorphism algebra 
${\rm End}^0(\{A_\chi\})$. The product $\prod_{f, \epsilon(f)=\chi} A_f$, where
$A_f$ is the abelian variety associated by Shimura to the cuspidal Hecke eigenform $f$ of weight $2$
and $f$ runs over all such normalized eigenforms with Nebentypus character $\chi$,
belongs to this isogeny class. Here we fix 
another representative $A_\chi$ of this isogeny class by picking for each prime $l$
the Galois stable lattice 
$
({\rm H}^1_{et}(A\otimes_\Q{\bar \Q}, \Z_l)\otimes_{\Z_l}{\Z_l}[\zeta_r]\chi^{-1})^\Gamma\subset 
V^\chi_l
$.
The abelian variety $A_\chi$ supports multiplication by $\Z[\zeta_r]$, i.e an embedding 
$\Z[\zeta_r]\subset {\rm End}(A_\chi)$. Let us denote by $\A_\chi$ the N\'eron model of $A_\chi$ over $\Z$. 
Then we can show that
${\rm Lie}(\A_\chi)\otimes_\Z\Z[\frac{1}{p-1}]\simeq \left(\overline{S_2(\Gamma_1(p), \Z[\zeta_r])_\chi}\right)^*\otimes_\Z\Z[\frac{1}{p-1}]$
where ${}^*$ denotes the $\Z$-dual of $\Z[\zeta_r]$-modules.
If $r$ is prime we have the conjectural identity
\begin{equation*}
\overline{\theta_2[\P_\chi]}=[\Sha(A_\chi)]+[\Phi_p(A_\chi)]-\overline{ [A_\chi(\Q)]}-[A_\chi(\Q)] 
\end{equation*}
in $\Gr_0(\Z[\zeta_r, \frac{1}{p-1}])/\{\hbox{\rm free modules}\}= {\rm Cl}(\Z[\zeta_r,\frac{1}{p-1}])$.
}
\end{Remarknumb}


\begin{thebibliography}{C-P-S1 }

\bibitem[A]{} A. Abb\`es, The Grothendieck-Ogg-Shafarevich formula for arithmetic surfaces,
J. Alg. Geom. 9 (2000), 529--576.



\bibitem[BFQ] {} P. Baum,  W. Fulton,  G. Quart: Lefschetz-Riemann-Roch for singular varieties. Acta Math. {\bf 143} (1979), no. 3-4, 193--211.



\bibitem[BLR] {} S. Bosch,  W. L\"utkebohmert, M. Raynaud: {\sl  N\'eron models.} Ergebnisse der Mathematik und ihrer Grenzgebiete (3), 21. Springer-Verlag, Berlin, 1990. 

\bibitem[Br] {} L. Breen: {\sl Fonctions th\^eta et th\'eor\`eme du cube.}
Lecture Notes in Mathematics, Vol. 980. Springer-Verlag, Berlin-New
 York, 1983.

\bibitem[BG] {} J. Buhler, B. Gross:
Arithmetic on elliptic curves with complex multiplication. II. 
Invent. Math. 79 (1985), no. 1, 11--29.

\bibitem[BF] {} D. Burns, M. Flach: Tamagawa numbers for motives with (non-commutative) coefficients. Doc. Math. 6 (2001), 501--570 (electronic).


\bibitem[CNT] {} Ph. Cassou-Nogu\`es, M. J. Taylor: Un \'el\'ement de Stickelberger quadratique.   J. Number Theory 37 (1991), no. 3, 307--342.
 
\bibitem[CW] {} C. Chevalley, A. Weil: \"Uber das Verhalten der Integrale 1.  Gattung bei Automorphismen des Funktionenk\"orpers, Abh. Math. Sem. Univ. Hamburg 10 (1934), 358--361; Zbl. 9, 160. 

\bibitem[C] {} T. Chinburg: Galois structure of de Rham cohomology of tame covers of
schemes, Annals of Math., {\bf 139} (1994), 443--490. Corrigendum, Annals of Math., {\bf 140} (1994), 251. 


\bibitem[CEPT1] {} T. Chinburg, B. Erez, G. Pappas, M. J. Taylor:
Tame actions of group schemes: integrals and slices, Duke Math. J. {\bf
82} (1996), 269--302.


\bibitem[CEPT2] {}  ------------------------------ :   $\varepsilon $%
-constants and the Galois structure of de Rham cohomology, Annals of Math.,
{\bf 146} (1997), 411--473. 



\bibitem[CEPT3] {} ----------------------------- :   Riemann-Roch type theorems for
arithmetic schemes with a finite
     group action, 
     J. Reine Angew. Math. {\bf 489}
(1997), 151--187.


\bibitem[CES] {} B. Conrad, B. Edixhoven, W. Stein:  $J\sb 1(p)$ has connected fibers. Doc. Math. 8 (2003), 331--408 (electronic). 

\bibitem[De] {} P.  Deligne: Le d\'{e}terminant de la cohomologie, in {\sl Currents Trends in 
Arithmetical Algebraic Geometry.} Contemporary Mathematics,
Vol. {\bf 67}, A.M.S. (1987).


\bibitem[DR] {} P. Deligne, M. Rapoport:
Les sch\'emas de modules de courbes elliptiques.  
Modular functions of one variable, II (Proc. Internat. 
Summer School, Univ. Antwerp, Antwerp, 1972), pp. 143--316. Lecture Notes in Math., Vol. {\bf 349}, 
Springer, Berlin, 1973. 



\bibitem[Du] {} F. Ducrot:  Cube structures and intersection bundles. J. Pure Appl. Algebra {\bf 195}, no. 1, (2004),
33--73.
  

\bibitem[F] {} M. Flach: 
The equivariant Tamagawa number conjecture: a survey. With an appendix by C. Greither. Contemp. Math., 358, Stark's conjectures: recent work and new directions, 79--125, Amer. Math. Soc., Providence, RI, 2004.  

 

\bibitem[Fr] {} A. Fr\"ohlich:  {\sl Galois Module Structure of Algebraic
Integers.} Ergebnisse der  Mathematik und 
ihrer Grenzgebiete.  (3) 1.  Springer-Verlag, Berlin, 1983.

\bibitem[Fu] {}  W. Fulton: 
{\sl Intersection theory.} $2^{\rm nd}$ edition. Ergebnisse der Mathematik und 
ihrer Grenzgebiete. 3. Folge. A Series of Modern Surveys in Mathematics, 2. Springer-Verlag, Berlin, 1998. 


 \bibitem[G] {} B. Gross: On the conjecture of Birch and Swinnerton-Dyer for elliptic curves with complex multiplication. Number theory related to Fermat's last theorem (Cambridge, Mass., 1981), pp. 219--236, Progr. Math., 26, Birkhäuser, Boston, Mass., 1982. 

\bibitem[GM] {} A. Grothendieck, J. P. Murre: 
The tame fundamental group of a formal neighbourhood of a divisor with normal crossings on a scheme. 
Lecture Notes in Mathematics, Vol. {\bf 208}. 
Springer-Verlag, Berlin-New York, 1971.  

\bibitem[GS] {} H. Gillet,  C. Soul\'e:
An arithmetic Riemann-Roch theorem. 
Invent. Math. {\bf 110} (1992),   473--543.


\bibitem[I]{} L. Illusie, Th\'eorie de Brauer et caract\'eristique d'Euler-Poincar\'e, Ast\'erisque
{\bf 82 - 83} (1978 - 79), 173--219.
 

\bibitem[KL] {} D. Kubert, S. Lang: {\sl  Modular units.} Grundlehren der Mathematischen Wissenschaften, 244. Springer-Verlag, New York-Berlin, 1981. 

\bibitem[KM] {} N. Katz, B. Mazur:
{\sl Arithmetic moduli of elliptic curves. }
Annals of Mathematics Studies, 108. 
Princeton University Press, Princeton, NJ, 1985. 

\bibitem[KS] {} N. Katz, P. Sarnak: Zeroes of zeta functions and symmetry.  Bull. Amer. Math. Soc. (N.S.)  36  (1999),  no. 1, 1--26. 

\bibitem[KMu] {} 
F.~Knudsen, D.~Mumford:
 The Projectivity of the Moduli space of Stable Curves (I):
  Preliminary on ``det" and ``Div".
Math. Scand., {\bf 39} (1975), 19--55.

\bibitem[La1]{} S. Lang:  {\sl Algebraic Number Theory.} second edition, Springer,
New York, 1994.

\bibitem[La2] {} ---------: 
{\sl Introduction to Arakelov theory. }
Springer-Verlag, New York, 1988. 

\bibitem[LLR] {} Q. Liu, D. Lorenzini, M. Raynaud: N\'eron models, Lie algebras, and reduction of curves of genus  
one. Invent. Math. {\bf 157}  (2004), 455--518.   



\bibitem[Ma] {}  B. Mazur: Modular curves and the Eisenstein ideal. Inst. Hautes 
\'Etudes Sci. Publ. Math. No. {\bf 47} (1977), 33--186. 



\bibitem[McC] {} L. McCulloh: Galois module structure of abelian extensions. 
J. Reine Angew. Math. {\bf 375/376} (1987), 259--306.



\bibitem[P1] {} G. Pappas: Galois modules and the Theorem of the Cube. Invent. Math.
{\bf 133} (1998), 193--225.

\bibitem[P2] {} ------------:  Cubic structures and ideal class groups.  
Ann. Sci. \'Ecole Norm. Sup. (4) {\bf  38} (2005), no. 3, 471--503.

\bibitem[P3] {} ------------: Galois module structure of unramified covers.
Preprint.

\bibitem[P4] {} ------------: Galois module structure and the $\gamma$-filtration.   
Compositio Math. {\bf 121} (2000), no. 1, 79--104.


\bibitem[P5] {} ------------: Chow groups and cubic structures. 
Preprint.

\bibitem[Ra] {} M. Raynaud:  Anneaux locaux hens\'eliens. 
 Lecture Notes in Mathematics, Vol. {\bf 169} (1970), Springer-Verlag, Berlin-New York. 
 

 
\bibitem[Ri] {}  D. S. Rim: Modules over finite groups. Annals of Math.,  {\bf 69} (1959), 700--712. 

\bibitem[Ro1] {} P. Roberts: Abelian extensions of regular local rings. 
Proc. Amer. Math. Soc. {\bf 78} (1980), no. 3, 307--310. 

 \bibitem[Ro2] {} -------------: {\sl Multiplicities and Chern classes in local algebra.} 
Cambridge Tracts in Mathematics, 133.
Cambridge University Press, Cambridge, 1998.



\bibitem[SGA4] {} {\sl Th\'eorie des topos et cohomologie \'etale des sch\'emas. } 
Dirig\'e par M. Artin, A. Grothendieck, et J.~L.~Verdier. Avec la collaboration de N.~Bourbaki, P.~Deligne
et B.~Saint-Donat. Lecture Notes in Mathematics, Vol. 269, 270, 305. Springer-Verlag, Berlin-New York, 1972-73.

\bibitem[SGA6] {} {\sl Th\'eorie des intersections et th\'eor\`eme de Riemann-Roch.}  
 Dirig\'e par P.~Berthelot, A.~Grothendieck et L.~Illusie. 
Avec la collaboration de D.~Ferrand, J.~P.~Jouanolou, O.~Jussila, S.~Kleiman, M.~Raynaud et 
J.~P.~Serre. Lecture Notes in Mathematics, Vol. 225. Springer-Verlag, Berlin-New York, 1971.


\bibitem[T1] {} R. W. Thomason: Lefschetz-Riemann-Roch theorem and coherent trace formula. Invent. Math. 85 (1986), no. 3, 515--543. 

\bibitem[T2] {} -----------------: Une formule de Lefschetz en $K$-th\'eorie \'equivariante alg\'ebrique.   Duke Math. J. 68 (1992), no. 3, 447--462.

\bibitem[V] {} I. Vidal, Th\'eorie de Brauer et conducteur de Swan.  J. Algebraic Geom. 13 (2004), no. 2, 349--391.

\bibitem[Wa] {} L. Washington:   {\sl Introduction to 
cyclotomic fields.} Graduate Texts in Mathematics, 83. Springer-Verlag, New York, 1982. xi+389 pp.
\end{thebibliography}
\end{document}